\documentclass[11pt,reqno]{amsart}
\usepackage{a4wide}
\usepackage{euscript,amsmath,amssymb,amsbsy,mathabx}
\usepackage[arrow,matrix,curve]{xy}
\usepackage{array,longtable,graphicx,enumitem,url}

\usepackage[colorlinks=true,linkcolor=black,urlcolor=black,citecolor=black]{hyperref}

\newcounter{msct}[section]\renewcommand{\themsct}{\thesection.\arabic{msct}}

\newenvironment{m-theorem}{\vskip5pt\refstepcounter{msct}\trivlist \itemindent 0pt \item[\hskip\labelsep\bf Theorem~\themsct] \it\ignorespaces}{\endtrivlist\vskip3pt}

\newenvironment{m-proposition}{\vskip5pt\refstepcounter{msct}\trivlist \itemindent0pt \item[\hskip\labelsep\bf Proposition~\themsct] \it\ignorespaces}{\endtrivlist\vskip3pt}

\newenvironment{m-corollary}{\vskip5pt\refstepcounter{msct}\trivlist \itemindent 0pt \item[\hskip\labelsep\bf Corollary~\themsct] \it\ignorespaces}{\endtrivlist\vskip3pt}

\newenvironment{m-lemma}{\vskip5pt\refstepcounter{msct}\trivlist \itemindent 0pt \item[\hskip\labelsep\bf Lemma~\themsct] \it\ignorespaces}{\endtrivlist\vskip3pt}

\newenvironment{m-definition}{\vskip5pt\refstepcounter{msct}\trivlist \itemindent0pt \item[\hskip\labelsep\bf Definition~\themsct] \ignorespaces}{\endtrivlist\vskip5pt}

\newenvironment{m-notation}{\vskip5pt\refstepcounter{msct}\trivlist \itemindent0pt \item[\hskip\labelsep\bf Notation~\themsct] \ignorespaces}{\endtrivlist\vskip5pt}

\newenvironment{m-example}{\vskip5pt\refstepcounter{msct}\trivlist \itemindent0pt \item[\hskip\labelsep\bf Example~\themsct] \ignorespaces}{\endtrivlist\vskip5pt}

\newenvironment{m-remark}{\vskip5pt\refstepcounter{msct}\trivlist \itemindent0pt \item[\hskip\labelsep\bf Remark~\themsct] \ignorespaces}{\endtrivlist\vskip5pt}

\newenvironment{m-question}{\vskip3pt\refstepcounter{msct}\trivlist \itemindent0pt \item[\hskip\labelsep\bf Question.] \ignorespaces}{\endtrivlist\vskip5pt}

\newenvironment{thm-nono}[1]{\vskip5pt\trivlist \itemindent 0pt \item[\hskip\labelsep\bf Theorem~{\rm\mbox{#1}}] \it\ignorespaces}{\endtrivlist\vskip5pt}
\newenvironment{lm-nono}[1]{\vskip5pt\trivlist \itemindent0pt \item[\hskip\labelsep\bf Lemma~{\rm\mbox{#1}}] \it\ignorespaces}{\endtrivlist\vskip5pt}
\newenvironment{conj-nono}[1]{\vskip5pt\trivlist \itemindent0pt \item[\hskip\labelsep\bf Conjecture~{\rm\mbox{#1}}] \it\ignorespaces}{\endtrivlist\vskip5pt}
\newenvironment{m-thank}{\vskip5pt\trivlist \itemindent0pt \item[\hskip\labelsep\it Acknowledgments] \ignorespaces}{\endtrivlist\vskip5pt}

\newenvironment{m-proof}{\vskip3pt\trivlist \itemindent0pt \item[\hskip\labelsep\it Proof.] \ignorespaces}{\hfill$\Box$\endtrivlist\vskip5pt}

\newenvironment{m-asmp}{\vskip5pt\trivlist \itemindent0pt \item[\hskip\labelsep\bf Assumption.] \ignorespaces}{\hfill\endtrivlist\vskip5pt}

\newcounter{meqn}[section]\renewcommand{\themeqn}{\thesection.\arabic{meqn}}
\newenvironment{m-eqn}[1]{\vskip5pt\refstepcounter{meqn}%
\trivlist\itemindent0pt\item[]\ignorespaces%
\hfill $\displaystyle #1$\hfill\hbox{\rm(\themeqn)}}{\endtrivlist\vskip5pt}

\newcommand{\bibauth}[2]{\textrm{{#1}~{#2}}}
\newcommand{\bibtitl}[1]{\textit{#1}.}
\newcommand{\bibjnyp}[4]{\textrm{#1} \textbf{#2} (#3), #4.}
\newcommand{\bibinbook}[4]{In: \textrm{#1}\textrm{, #2}\textrm{, #3}\textrm{, #4}.}
\newcommand{\bibbook}[4]{\textit{#1}. {#2} {#3}, {#4}.}

\let\rar\rightarrow
\let\lar\longrightarrow
\let\hra\hookrightarrow
\let\mt\mapsto
\let\lmt\longmapsto
\let\dashto\dashrightarrow
\newcommand\arr[1]{{\vrule height3.25pt depth-2.6pt width #1\kern-2ex\hbox{$\rightarrow$}}}

\let\euf\EuScript 
\let\cal\mathcal
\let\mbb\mathbb

\let\mfrak\mathfrak
 
\DeclareFontFamily{OT1}{rsfs}{}
\DeclareFontShape{OT1}{rsfs}{n}{it}{<->rsfs10}{}
\DeclareMathAlphabet{\crl}{OT1}{rsfs}{n}{it}

\numberwithin{equation}{section}\numberwithin{figure}{section} 

\makeatletter
\renewcommand{\p@enumi}{}
\renewcommand{\p@enumii}{}
\makeatother

\newcommand\bone{{1\kern-0.57ex\rm l}}
\newcommand{\bbA}{{\mbb A}}
\newcommand{\cA}{{\cal A}}
\newcommand{\eA}{{\euf A}}
\newcommand{\eB}{{\euf B}}
\newcommand{\cD}{{\cal D}}

\newcommand\beF{{{\euf F}}}

\newcommand\eG{{\euf G}}
\newcommand\beG{{{\euf G}}}

\newcommand{\Gl}{{\mathop{\rm GL}\nolimits}}
\newcommand\eI{{\euf I}}
\newcommand{\eK}{{\euf K}}
\newcommand\eO{{\euf O}}
\newcommand{\cQ}{{\cal Q}}
\newcommand{\eQ}{{\euf Q}}
\newcommand{\eS}{{\euf S}}
\newcommand{\eT}{{\euf T}}
\newcommand{\del}{{\partial}}

\newcommand{\ev}{\mathop{{\rm ev}}\nolimits}
\let\veps\varepsilon 
\newcommand\Aut{\mathop{\rm{Aut}}\nolimits}
\newcommand\dne{\mathop{\rm{end}}\nolimits}
\newcommand\End{\mathop{\rm{End}}\nolimits}
\newcommand\cEnd{\operatorname{{\cal E}\kern-1pt \textit{nd}\kern1pt}}
\newcommand\ext{\mathop{\rm{ext}}\nolimits}
\newcommand\Ext{\mathop{\rm{Ext}}\nolimits}
\newcommand\cExt{\operatorname{{\cal E}\kern-1pt\textit{xt}\kern1pt}}
\newcommand{\ges}{{\geqslant}}
\newcommand{\cH}{{\cal H}}
\newcommand{\iph}{{\cal I}}
\newcommand{\Hom}{\mathop{{\rm Hom}}\nolimits}
\newcommand{\cHom}{\mathop{{\cal H}om}\nolimits}
\newcommand{\Img}{\mathop{\rm Im}\nolimits}
\newcommand{\invq}{{/\!/}}
\newcommand{\IVB}{{M\!I}_\ppp}
\newcommand{\kz}{{{\kappa}\!z}}
\newcommand{\bi}{{\boldsymbol i}}\newcommand{\bj}{{\boldsymbol j}}\newcommand{\bk}{{\boldsymbol k}}
\newcommand{\Ker}{{\mathop{\rm Ker}\nolimits}}
\newcommand{\bbL}{{\mbb L}}
\newcommand{\LL}{{I\kern-.8ex L}}
\let\lda\lambda
\let\les\leqslant 
\newcommand{\lcit}{{\textit{loc.\,cit.}}}

\newcommand{\lft}{{\rm left}}\newcommand{\rgt}{{\rm right}}
\let\nit\noindent 
\newcommand{\mi}{{M\!I}}
\newcommand{\mm}{{M}}
\newcommand{\pl}{{{\rm p}_{l}}}
\newcommand{\pr}{{{\rm p}_{r}}}
\newcommand{\ppp}{{\mbb P^3}}\newcommand{\pp}{{\mbb P^2}}
\newcommand{\PGL}{{\mathop{\rm PGL}\nolimits}}

\newcommand{\rR}{{\mathrm R}}
\newcommand{\res}{{\mathop{\rm rs}\nolimits}}
\newcommand\ort{\mathrel{{\vrule width4pt height0.4pt depth0pt\vrule width0.4pt height6pt depth0pt\,}}}
\newcommand\ouset[3]{{\overset{#2}{\underset{#1}#3}}}
\newcommand{\rst}{{\!\upharpoonright}}
\let\si\sigma 
\let\surj\twoheadrightarrow
\let\srel\stackrel 
\newcommand{\svb}{{\rm svb}}
\newcommand{\bbT}{{\mbb T}}
\newcommand{\cU}{{\cal U}}
\newcommand{\cV}{{\cal V}}
\newcommand{\ux}{{\underline{x}}}
\newcommand{\uz}{{\underline{z}}}
\let\tld\tilde 
\let\wtld\widetilde

\title[Properties of mathematical instantons on $\mbb P^3$]{Rationality and categorical properties of the moduli of instanton bundles on the projective 3-space}
\author[M.~Halic, R.~Tajarod]{Mihai~Halic, Roshan~Tajarod}
\email{mihai.halic@gmail.com, roshan.tajarod@gmail.com}
\begin{document}
\begin{abstract}
We prove the rationality and irreducibility of the moduli space of  mathematical instanton vector bundles on $\ppp$, of arbitrary rank and charge. In particular, the result applies to the rank-$2$ case. This problem was first studied by Barth, Ellingsrud-Str{\o}mme, Hartshorne, Hirschowitz-Narasimhan in the late 1970s. We also show that the mathematical instantons of variable rank and charge form a monoidal category. The proof is based on a careful analysis of the Barth-Hulek monad-construction and on a detailed description of the moduli space of (framed and unframed) stable bundles on Hirzebruch surfaces.
\end{abstract}
\maketitle 

\section*{Introduction}
The interest in rank-$2$ instanton bundles on the three-dimensional projective space, with Chern classes $c_1=0,c_2=n$, has its origins in the articles of Atiyah \textit{et al.}~\cite{adhm,ahs,atiyah-ward}, Barth-Hulek~\cite{brth+hulk}, Drinfel'd-Manin~\cite{dr-ma}, and Hartshorne~\cite{hart-vb-P3, hart-inst, hart-probl}, which in turn were motivated by work of 't~\!Hooft~\cite{hooft} and Polyakov~\cite{plkv}. The geometry of their moduli spaces, such as smoothness and irreducibility, has been intensively investigated: see~\cite{brth-c=4,dons,sin-tra,ctt,jmt,tikh-maincomp} for the rank-$2$ case and~\cite{atiy-2+4,bruz-mark-tikh,bruz-mark-tikh2,dons,fre-jar,jar-ver-twis, jar-ver-tri} for generalizations to higher rank bundles on the projective space $\ppp$. 

The issue regarding the rationality of the moduli spaces of instanton bundles was raised by Hartshorne~\cite{hart-vb-P3}, and turned out to be difficult. So far, it's known that these varieties are rational for rank-$2$ and charge $n=2,4,5$, due to works of Hirschowitz-Narasimhan~\cite{hirs+nars}, Ellingsrud-Str{\o}mme~\cite{ellin+strm}, and Katsylo~\cite{kats}; for $n\ges 6$, the issue remained open, in spite of efforts~\cite{tikh-odd,tikh-even}. Beorchia-Franco~\cite{beo+fra} proved that the moduli space of 't~\!Hooft instantons --those which possess a section at the first twist-- is irreducible and rational. Let us emphasize that the techniques are specific to the rank-two case. 

Our goal is to address the rationality issue mentioned above in \textit{arbitrary rank}. Most of the literature is dedicated to the rank-$2$ case, but non-abelian gauge theories  for $SU(r)$ are frequent in physics. The ADHM-construction and the Penrose-transform relating Hermitian vector bundles with self-dual connections over the sphere $\mbb S^4$, on the one hand, to certain holomorphic vector bundles on $\ppp$, on the other hand, apply in this general setting~\cite{adhm,ahs}. Therefore we believe that our \emph{unified treatment} of the \emph{arbitrary rank case} yields additional interest to this article. 

Let us precise that our definition of mathematical instantons assumes the \textit{existence of a trivializing line}; rank-$2$, semi-stable vector bundles with $c_1=0$ automatically satisfy this condition. The first, possibly unexpected, result is that the instanton-property of vector bundles (of variable rank and charge) --that is, vanishing of the $H^1$- and $H^2$-cohomology groups of their $(-2)$-twist-- is closed under tensor product. 
This property is classical for the Yang-Mills (YM) instantons, it follows from the Penrose correspondence~\cite{ahs}. Our work generalizes this property to arbitrary mathematical instantons. 

\begin{thm-nono}{\textbf{A (cf.~\ref{thm:categ}})} 
The category $(\IVB,\otimes)$, whose objects are mathematical instanton bundles on $\ppp$, is a self-dual monoidal category. The subcategory $\IVB^{\rm ps}$ consisting of poly-stable bundles (finite direct sums of stable objects) possesses the structure of a multi-tensor category. In particular, Schur powers (representations) preserve mathematical instantons. 
\end{thm-nono}

To our knowledge, it hasn't been observed so far that the moduli spaces of instanton bundles aggregate into an algebraic object (see~\cite{egno} for definitions). However, it's a pleasant surprise that tensor products of physical instantons --those arising from the ADHM-construction-- are relevant in QCD~\cite{cgt,dri-man,jack,be-st,tem}. We very briefly discuss this matter in the last section.

The invariance of the instanton-property for tensor products has two important down-to-earth consequences: the moduli space of mathematical instanton bundles of given rank $r$ and charge $n$ has expected dimension, given by the Riemann-Roch formula, and the general instanton on $\ppp$ is uniquely determined by its restriction to either a wedge of two $2$-planes or to a smooth quadric. This restriction property has remarkable implications. 

First, it reduces the initial problem to the study of semi-stable vector bundles on Hirzebruch surfaces. Explicit computations become possible and we obtain the following result.

\begin{thm-nono}{\textbf{B (cf.~\ref{thm:rtl-frame})}}
The moduli space of (unframed and framed) semi-stable vector bundles on a Hirzebruch surface, with Chern classes $c_1=0, c_2=n$ is rational. 
\end{thm-nono}
We emphasize the \emph{novelty} even in this \emph{$2$-dimensional setting}: the moduli spaces of stable bundles on Hirzebruch surfaces are known to be unirational for all $n$ (cf.~\cite{lq}) and rational for $n\gg0$ (but without explicit bounds, cf.~\cite{c-mr}). 
The result allows to conclude:  

\begin{thm-nono}{\textbf{C (cf.~\ref{thm:irred})}}
The moduli space $\mi_\ppp(r;n)_{(\lda)}$ of mathematical ($\lda$-framed) instanton bundles on $\ppp$, of rank $r$, with Chern classes $0,n,0$, is irreducible of dimension $4rn-r^2+1$ (resp. $4rn$), and is {\emph{rational}}. In particular, these properties hold in the rank-$2$ case.
\end{thm-nono}

Second, the restriction property implies an unexpected relation between the mathematical- and the YM-instantons. Let $\iph_{\mbb{CP}^3}(r;n)_{\rm line}$ be the moduli space of $SU(r)$ YM-instantons of charge $n$, with framing along a real line in $\mbb{CP}^3$; Atiyah~\cite{atiy-2+4} endowed it with a complex-analytic structure. We give an alternative, algebraic proof to a question raised by Atiyah~\cite[Problem 22]{hart-probl}, solved by Donaldson~\cite{dons}, and we deduce a surprising consequence: roughly speaking, mathematical instantons of charge $n$ are the same as Yang-Mills instantons of same rank and charge $2n$. In terms of the monad/ADHM-construction, this is completely unclear.

\begin{thm-nono}{\textbf{D (cf.~\ref{thm:ym},~\ref{thm:n2n})}} 
\begin{enumerate}[leftmargin=5ex]
\item 
$\iph_{\mbb{CP}^3}(r;n)_\lda\cong\mm_{\mbb{CP}^2}(r;n)_\lda$ is a rational, complex quasi-projective variety of dimension $2rn$. 
\item 
Let $\lda,\lda'\subset\mbb{CP}^3$ be two intersecting lines. 
There is an (algebraic) open immersion $\mi_{\mbb{CP}^3}(r;n)_{\lda\cup\lda'}\to\iph_{\mbb{CP}^3}(r;2n)_{\rm line}$ commuting with direct sums, tensor products.
\end{enumerate}
\end{thm-nono}

We conclude this introduction with a survey of the article. The techniques are cohomo\-logi\-cal in nature but, compared to existing approaches, we dissect homomorphisms between relevant cohomology groups to understand their action.
\begin{itemize}[leftmargin=3ex]
\item 
The analysis of the properties of instantons requires an \emph{in-depth understanding} of the Barth-Hulek construction~\cite{brth+hulk}. We show that their monad is determined by a `universal' diagram which is based on Beilinson's resolution of the diagonal in $\ppp\times\ppp$. 
The matter is essential to explicitly determine the various homomorphisms between cohomology groups induced by the monad, and allows performing computations. The main outcome is the invariance of the instanton-property under tensor products (Theorem~A).  
\item 
We analyse the geometry of the moduli space of instantons on $\ppp$ by restricting them to a wedge of planes, resp. a smooth quadric. The invariance property implies that these \textit{restriction maps are birational} (cf. Theorem~\ref{thm:birtl}). 

To appreciate the strength of this statement, note that~\cite{lq, c-mr} immediately imply the unirationality (for all $n$) and the rationality (for $n\gg0$) of $\mi_\ppp(r;n)$. Even the former, weaker property has been out of the reach of current approaches.
\item 
In Section~\ref{sct:hirz},  we apply methods developed by the first named author~\cite{hlc-sheaves} to prove the rationality of the moduli spaces of stable bundles on Hirzebruch surfaces (Theorem~B). 
\item 
We deduce Theorem~D, relating mathematical and Yang-Mills instantons, in Section~\ref{sct:roots}. 
\end{itemize}
\nit We work over an algebraically closed field $\Bbbk$ of characteristic zero. For shorthand, the symbol `$\hra$' indicates monomorphisms and `$\surj$' epimorphisms;  short exact sequences are denoted $A\hra C\surj B$. The notation `$\ev$' will stand for evaluation maps, `$\ort$' for contraction (pairing), `$\rst, \res$' for restrictions, and `$\del$' for boundary maps in cohomology.


\section{The framework}\label{sct:frame}

\begin{m-definition}\label{def:FonP3} 
A \emph{mathematical instanton} on $\mbb P^3$ of rank $r$ and charge $n$ is a vector bundle $\beF$ which satisfies the following properties:
\begin{enumerate}[leftmargin=5ex]
\item[(i)] 
${\rm rank}(\beF)=r$, $c_1(\beF) = 0, c_2(\beF)=n, c_3(\beF)=0$;
\item[(ii)] 
its \emph{restriction} to a (general) line $\lda_{gen}\subset\mbb P^3$ is \emph{trivializable} (so $\beF$ is slope semi-stable);
\item[(iii)] 
it satisfies the \emph{instanton condition}: $H^1(\beF(-2)) = H^2(\beF(-2)) = 0.$
\end{enumerate}
\end{m-definition}

\begin{m-remark}\label{rmk:FonP3}
\begin{enumerate}[leftmargin=5ex]
\item 
For $r=2$, semi-stable vector bundles with $c_1=0$ automatically satisfy the restriction property~(ii), by the Grauert-M\"ullich theorem.
\item 
For arbitrary $r$, Barth-Hulek~\cite{brth+hulk} showed that vector bundles $\beF$ satisfying (i)-(iii) are the cohomology of a monad, which is determined up to isomorphism: 
\begin{m-eqn}{
\underbrace{H^1(\beF^\vee(-1))^\vee}_{\cong H^2(\beF(-3))\cong\Bbbk^n}\otimes\eO_\ppp(-1)
\to\eO_\ppp^{\oplus r+2n}\to\underbrace{H^1(\beF(-1))}_{\cong\Bbbk^n}\otimes\eO_\ppp(1).
}\label{eq:bh-mnd3}
\end{m-eqn}
For our approach it's important to pinpoint a \emph{unique} monad which yields $\beF$. 
\item 
The Riemann-Roch yields $h^1(\beF)=2n-r$, so one has $2n\ges r$. We assume $n\ges r$ throughout. The reason for this hypothesis is given in Lemma~\ref{lm:stable}.
\end{enumerate}
\end{m-remark}

\begin{m-notation}\label{not:moduli} 
We consider the following quasi-projective varieties:
\begin{enumerate}[leftmargin=5ex]
\item 
$\mi_{\ppp}(r;n),$ the moduli space of instanton vector bundles of rank $r$ and charge $n$. It is a non-empty open subset of the moduli space of slope semi-stable sheaves on $\ppp$.

Let $\beF$ be an instanton and $\lda\subset\ppp$ a line, such that the restriction $\beF_\lda$ is trivializable. A \emph{framing} of $\beF$ along $\lda$ is an isomorphism $\alpha_\lda:\beF_\lda\to\eO_\lda^{\oplus r}$. The frames $\alpha_\lda, \alpha'_\lda$ of $\beF,\beF'$, respectively, are equivalent if there is a commutative the diagram as below:
$$
\xymatrix@R=1.5em@C=3em{
\beF\ar[d]_-{a}^-{\cong}\ar[r]^-{\res}
&\beF_\lda\ar[d]_-{a_\lda}^-{\cong}\ar[r]^-{\alpha_\lda}|-{\;\cong\;}
&\eO_\lda^{\oplus r}\ar@{=}[d]
\\ 
\beF'\ar[r]^-{\res}
&\beF'_\lda\ar[r]^-{\alpha'_\lda}|-{\;\cong\;}
&\eO_\lda^{\oplus r}
}
$$
Thus two frames of a stable instanton $\beF$ are equivalent if they differ by a multiplicative factor. Let $\mi_\ppp(r;n)_\lda$ be the moduli space of framed vector bundles.

\smallskip\item 
$\mm_{\pp}(r;n)$, resp. $\mm_{\pp}(r;n)_\lda$, the moduli space of rank-$r$, slope semi-stable, resp. framed along the line $\lda$, vector bundles on $\pp$, with $c_1=0, c_2=n$; for simplicity, we call them \emph{$\pp$-instantons}. In the framed case, the vector bundles are trivializable along $\lda$. The moduli spaces are irreducible and $(2rn-r^2+1)$-, resp. $2rn$-dimensional, see Section~\ref{sct:hirz}.
\smallskip\item 
For $2$-planes $\cD, \cH\subset\ppp$ and $\lda:=\cD\cap \cH$, let  
$$
\mm_{\cD\cup \cH}(r;n)_\lda:=\mm_\cD(r;n)_\lda\times \mm_\cH(r;n)_\lda,
$$
be the variety of semi-stable, framed vector bundles on $\cD\cup \cH$; the frames of the factors are used for gluing along $\lda$. Let 
$$
\mm_{\cD\cup \cH}(r;n):=\frac{\mm_{\cD\cup \cH}(r;n)_\lda}{\PGL(r)} =\frac{\mm_\cD(r;n)_\lda\times \mm_\cH(r;n)_\lda}{\PGL(r)},
$$
be the variety of semi-stable vector bundles on $\cD\cup \cH$, where the group acts diagonally on the frames. An element of $\mm_{\cD\cup \cH}(r;n)$ is stable if its restrictions to $\cD, \cH$ are so.

For a quadric $\cQ\cong\mbb P^1\times\mbb P^1$, let $\mm_\cQ(r;n)$ be the moduli space of vector bundles on $\cQ$, with $c_1=0, c_2=n$, which are semi-stable for $\eO_{\mbb P^1}(1)\boxtimes\eO_{\mbb P^1}(c),\,c>r(r-1)n$.
\end{enumerate}
\end{m-notation}

\begin{m-proposition}
The moduli space $\mi_\ppp(r;n)$ is non-empty. 
\end{m-proposition}

\begin{m-proof}
For $r=2$, consider the union $Z$ of $n+1$ disjoint lines in $\ppp$. The rank-$2$ vector bundle given by the Hartshorne-Serre construction along $Z$ fits into the exact sequence 
$$
\eO_\ppp(-1)\hra\beF_2^{HS}\surj\cal I_Z(1).
$$ 
(The construction appears in~\cite[Example 3.1.1]{hart-vb-P3}.) For $r>2$,  let we consider the rank-$r$ bundle 
$\euf F_r^{HS}:=\euf O_\ppp^{\oplus (r-2)}\oplus\,\euf F_2^{HS}.$ 
One easily verifies that $\beF_r^{HS}\otimes\eO_Z\cong\eO_Z^{\oplus r}$ and both $\beF_r^{HS}$,  $\cEnd(\beF_r^{HS})$ satisfy the conditions~(i)-(iii) in Definition~\ref{def:FonP3}. 
\end{m-proof}


\section{The monad construction revisited} \label{sct:monad}

Let $\beF$ be a mathematical instanton on $\ppp$.  Barth-Hulek~\cite[\S7]{brth+hulk} proved that it's isomorphic to the cohomology $\Ker(q)/\Img(\veps)$ of a (linear) monad 
\begin{align*}
\eO_\ppp(-1)^{n}\srel{\veps}\to\eO_\ppp^{r+2n}\srel{q}{\to}\eO_\ppp(1)^{n}.
\tag{$\varstar$}\label{eq:LM}
\end{align*}
For stable bundles, this is uniquely defined up to the natural $\tilde A := \Gl(n)\times\Gl(r+2n)\times\Gl(n)$ action on the terms. The diagonally embedded $\Bbbk^*$ acts trivially, so we get an $A:=\tilde A/\Bbbk^*$ action on the moduli space of monads~\eqref{eq:LM}. The latter is an open subset of the affine variety
\begin{m-eqn}{
Cplx_\ppp(r;n):=\{(\veps,q)\mid q\circ\veps=0\}\subset 
(Mat_{r+2n,n}\times Mat_{n,r+2n})\otimes H^0(\eO_\ppp(1)).
}\label{eq:cplx}
\end{m-eqn}
The dimension of its general irreducible component(s) is: 
$$
2\cdot n(r+2n)\cdot h^0(\eO_\ppp(1))-n^2\cdot h^0(\eO_\ppp(2))
=8n(r+2n)-10n^2=8rn+6n^2.
$$
Let $Mond_\ppp(r;n)$ be the open subset, of this dimension, corresponding to monads --$\veps$ is injective, $q$ is surjective-- such that their cohomology is trivializable along a general line. One gets an $A$-invariant (quotient) map 
$\,Mond_\ppp(r;n)\to\mi_\ppp(r;n).$ 
A dimension counting yields: 
$$
\dim Mond_\ppp(r;n) - \dim A = 8rn+6n^2 - (2n^2+(r+2n)^2) +1 = 4rn -r^2 +1.
$$
This is consistent with the (subsequent) cohomological calculation in Proposition~\ref{prop:EgI}. 

The main goal of this section is, for an instanton $\beF$, to determine a {\em canonically} associated monad whose cohomology is precisely $\beF$. By this, we mean that the isomorphism-ambiguity in~\eqref{eq:LM} should be at most $\Bbbk^*$, the automorphisms of stable bundles. 
Barth-Hulek's construction implies that the display of the monad whose cohomology is $\beF$ is: 
\begin{align*}
\xymatrix@R=1.5em@C=1.25em{
V_\beF\otimes\eO_\ppp(-1)\ar@{^(->}[r]^-{\veps}\ar@{=}[d]&
\;\;\eK_\beF\ar@{->>}[r]\ar@{^(->}[d]^{a^\vee}&
\;\beF\ar@{^(->}[d]&
\\ 
V_\beF\otimes\eO_\ppp(-1)\ar@{^(->}[r]^-{\wtld\veps}&
C_\beF\otimes\eO_\ppp\ar@{->>}[r]^-{b}\ar@{->>}[d]^-{\wtld q}&
\eQ_\beF\ar@{->>}[d]^-{q}&
\\ 
&W_\beF\otimes\eO_\ppp(1)\ar@{=}[r]&W_\beF\otimes\eO_\ppp(1).&
}
\tag{BH}\label{eq:VW}
\end{align*} 
The uniqueness part of the construction implies that two monads which determine the same vector bundle are in the same $\Gl(n)\times\Gl(n)$-orbit, where  
$$
\Gl(n)\times\Gl(n)=\frac{\Gl(n)\times\Bbbk^*\times\Gl(n)}{\Bbbk^*}\subset A
$$
acts on the extremities of~\eqref{eq:LM}. Let us enumerate several properties of~\eqref{eq:VW}. 

\begin{enumerate}[leftmargin=5ex]
\item 
There are \emph{canonical} identifications: 
\item[] 
$W_\beF=H^1(\beF(-1));$\quad  $V_\beF=H^2(\beF(-3))\srel{Serre}{\cong}H^1(\beF^\vee(-1))^\vee=(W_{\beF^\vee})^\vee.$ 
\item[] 
$C_\beF =H^0(\eQ_\beF)\cong H^0(\eK_\beF^\vee)^\vee.$

\item[] 
So $b$ is $H^0(\eQ_\beF)\otimes\eO_\ppp\to\eQ_\beF$ and $a^\vee$ is the dual of the evaluation map to $\eK_\beF^\vee$.

\item 
The diagram is obtained as follows:
\begin{enumerate}
\item 
The right-hand column is the extension defined by $\bone\in\End(W_\beF)=\Ext^1(W_\beF(1),\beF).$ 
\item 
The top row is defined by $\bone\in\End(V_\beF)=\Ext^1(\beF,V_\beF(-1)).$
\end{enumerate}
\item[] 
All these `rigidifications' still leave open the issue that~\eqref{eq:VW} is determined up to a $\Gl(n)\times\Gl(n)$-action, simply because all yield the same $\beF$. So, to achieve our goal, we need to explicitly determine homomorphisms $\veps_\beF, q_\beF$, depending naturally on $\beF$, which fit into the display; if an element of  $\Gl(n)\times\Gl(n)$ fixes the homomorphisms, too, then it's necessarily the identity. Explicit expressions for $\veps_\beF, q_\beF$ are necessary to understand what are the homomorphisms between various cohomology groups. 
\item  
The restriction $\beF_\cH=\beF\otimes\eO_\cH$ to a general hyperplane $\cH\cong\pp$ in $\ppp$ is still semi-stable, with the same Chern classes. The middle terms of the exact sequence  
$$
0=H^1(\beF(-2))\to H^1(\beF(-1))\srel{\res_\cH}{\lar} H^1(\beF_\cH(-1))\to H^2(\beF(-2))=0
$$
 are isomorphic. Thus the monad for $\beF_\cH,$ 
$$
\underbrace{H^1(\beF_\cH^\vee(-1))^\vee}_{\cong H^1(\beF_\cH(-2))}\!\otimes \eO_{\cH}(-1)\to\eO_{\cH}^{\oplus r+2n}\!\to\! H^1(\beF_\cH(-1))\!\otimes\!\eO_{\cH}(1),
$$  
is the restriction of~\eqref{eq:VW} to $\cH$. This reduction from $3$- to $2$-dimensions is essential for understanding the geometry of the moduli space of instanton bundles on $\ppp$. 
\end{enumerate}


\subsection{Koszul resolution}\label{ssct:kz}

Let $\lda$ be a trivialising line for $\beF$; denote $\eI_\lda$ its sheaf of ideals. By applying $\Hom(\cdot,\beF(-3))$ to the Koszul resolution $\eO(-2)\hra\eO(-1)^{\oplus 2}\surj\eI_\lda\,$, we obtain 
$$ 
0\to\underbrace{\Ext^1(\eO(-2),\beF(-3))}_{=H^1(\beF(-1))=W_\beF}
\;\ouset{\cong}{\kz_{\beF}}{\lar}\;
\underbrace{\Ext^2(\eI_\lda,\beF(-3))}_{=H^2(\beF(-3))= V_\beF}\to 0.
$$
On the right-hand side, we used $\beF_\lda\cong\eO_\lda^{\oplus r}$. Thus $\kz_\beF$ is the Yoneda-product (pairing) with the element of $\Ext^1(\eI_\lda,\eO(-2))$ defining the resolution.

We give now an alternative description of $\kz_\beF$ which yields a `formula' for its inverse. Let $\cD, \cH\subset\ppp$ be two hyperplanes intersecting along $\lda$. They determine the diagrams below; the Koszul resolution is the middle row of the second one: 
\begin{m-eqn}{
\begin{array}{c|c|c}
\scalebox{.9}{$
\xymatrix@R=1.5em@C=1em{
\;\;\eI_\cH\ar@{=}[r]\ar@{^(->}[d]&\;\;\eI_\cH\ar@{^(->}[d]
\\ 
\eI_\lda\ar@{^(->}[r]\ar@{->>}[d]&\eO_\ppp\ar@{->>}[d]
\\ 
\eO_\cH(-\lda)\ar@{^(->}[r] &\eO_\cH
}
$}
&
\scalebox{.9}{$
\xymatrix@R=1.5em@C=1em{
&\;\;\eI_\cH\ar@{=}[r]\ar@{^(->}[d]&\;\;\eI_\cH\ar@{^(->}[d]
\\ 
\eI_{\cD\cup \cH}\ar@{^(->}[r]\ar@{=}[d]&\eI_\cD\oplus\eI_\cH\ar@{->>}[d]\ar@{->>}[r]
&\eI_\lda\ar@{->>}[d]
\\ 
\eI_{\cD\cup \cH}\ar@{^(->}[r]&\eI_\cD\ar@{->>}[r]&\eO_\cH(-\lda)
}
$}
&
\scalebox{.9}{$
\xymatrix@R=1.5em@C=1em{
&\eI_{\cD\cup \cH}\ar@{^(->}[d]\ar@{^(->}[r] &\eI_\cH\ar@{^(->}[d]\ar@{->>}[r]&\;\eO_{\cD}(-\lda)
\\ 
&\eO_\ppp\ar@{=}[r]\ar@{->>}[d]
&\eO_\ppp\ar@{->>}[d]&
\\ 
\eO_\cD(-\lda)\ar@{^(->}[r]&\eO_{\cD\cup \cH}\ar@{->>}[r]&\eO_\cH&
}
$}
\end{array}
}\label{eq:IDH}
\end{m-eqn}
By taking the tensor product of the first two with $\beF(-1)$, we get the commutative diagram:
\begin{m-eqn}{
\xymatrix@R=1.5em@C=1.25em{
H^1(\beF(-1))\ar@/^4ex/[rr]|-{\;\kz_\beF\;}\ar[d]_-\cong^-{\res_\cH}
&H^1(\eI_\lda\otimes\beF(-1))\ar[l]^-{\cong}\ar[d]^-\cong\ar[r]_-\cong
&H^2(\beF(-3))\ar[d]_\cong^(.45){\del^{-1}}
\\ 
H^1(\beF_\cH(-1))&H^1(\beF_\cH(-2))\ar[l]^\cong\ar@{=}[r]&
H^1(\beF_\cH(-2)).
}
}\label{eq:kz-1}
\end{m-eqn}
The inverse of $\kz_\beF$ is obtained by following the lower edges: the restrictions to $\cH$ are isomorphisms and  $\kz_\beF^{-1}$ becomes the inclusion map $H^1(\beF_\cH(-2))\srel{\cong}{\to}H^1(\beF_\cH(-1))$. This observation will be essential later on.


\subsection{Beilinson resolution}\label{ssct:beil}

Let $\pl,\pr:\ppp\times\ppp\to\ppp$ be the projections onto the first and second factors,  and $\Delta\subset\ppp\times\ppp$ be the diagonal. 
The Euler sequence on $\ppp$ is: 
$$
\eO_\ppp(-1)\hra H^0(\eT_\ppp(-1))\otimes\eO_\ppp\srel{\ev}{\surj}\eT_\ppp(-1),\quad 
\Omega_\ppp^1(1)\hra H^0(\eO_\ppp(1))\otimes \eO_\ppp\srel{\ev}{\surj}\eO_\ppp(1).
$$
Note that $\pl_*(\eI_\Delta\otimes\pr^*\eO_\ppp(1))\cong\Omega^1_\ppp(1),$ and that $H^0(\eT_\ppp(-1)), H^0(\eO_\ppp(1))$ are dual to each other. 
Let 
$s\in H^0(\eT_\ppp(-1)\boxtimes\eO_\ppp(1))=\End H^0(\eO_\ppp(1))$ 
be the identity element. It transversally vanishes along $\Delta$, so one obtains the resolution of $\eI_\Delta$: 
$$\scalebox{.95}{$
0\to\ouset{}{3}{{\mbox{$\bigwedge$}}}
\bigl(\Omega_\ppp^1(1)\boxtimes\eO_\ppp(-1) \bigr)
\srel{s\ort}{\lar} 
\ouset{}{2}{{\mbox{$\bigwedge$}}}
\bigl(\Omega_\ppp^1(1)\boxtimes\eO_\ppp(-1) \bigr)
\srel{s\ort}{\lar} 
\kern-3ex\overbrace{\Omega_\ppp^1(1)\boxtimes\eO_\ppp(-1)}^{\pl_*(\eI_\Delta\otimes\pr^*\eO_\ppp(1))\boxtimes\eO_\ppp(-1)}\kern-3ex
\ouset{\ev}{s\ort}{\lar}\eI_\ppp\to0,
$}$$
where the homomorphisms are contractions with $s$. We consider the sheaf  
\begin{m-eqn}{
\eS:=\Ker(\ev)\cong\frac{\Omega^2_\ppp(2)\boxtimes\eO_\ppp(-2)}{\eO_\ppp(-1)\boxtimes\eO_\ppp(-3)}
\cong\frac{\eT_\ppp(-2)\boxtimes\eO_\ppp(-2)}{\eO_\ppp(-1)\boxtimes\eO_\ppp(-3)},
}\label{eq:eS}
\end{m-eqn}
which fits into the exact sequence $\;\eS\hra \Omega_\ppp^1(1) \boxtimes\eO_\ppp(-1)\surj\eI_\Delta.$ 

\begin{m-proposition}\label{prop:hom}
Let $\beF$ be an instanton vector bundle. The display~\eqref{eq:VW} is obtained by applying suitable derived $\cHom$-functors to the `universal' diagram below (\underbar{\emph{independent}} of $\beF$): 
\vskip0.5em\begin{m-eqn}{
\scalebox{.8}{$
\xymatrix@R=1em@C=1.5em{
\pr^*\eO_\ppp(1)\ar@/^4ex/@{->>}[rr] &\,\eI_\Delta\otimes\pr^*\eO_\ppp(1)\ar@{_(->}[l]&\eO_\Delta(1) 
&\genfrac{}{}{0pt}{1}{\text{apply}}{\rR\pr_*\rR\cHom(\,\cdot\,,\pl^*\beF(-3))} 
\\ 
\eS\otimes\pr^*\eO_\ppp(1)\,\ar@{^(->}[r]
&\pl^*\Omega^1_\ppp(1)\ar@{->>}[r]\ar@{->>}[u]
&\eI_\Delta\otimes\pr^*\eO_\ppp(1)\ar@{^(->}[d] 
&\genfrac{}{}{0pt}{1}{\text{apply}}{\rR\pr_*\rR\cHom(\pl^*\beF^\vee(1),\,\cdot\,)}
\\ 
&\eS\otimes\pr^*\eO_\ppp(1)\ar@{^(->}[u]
&\pr^*\eO_\ppp(1)\ar@/_10ex/@{->>}[uu] 
& 
\\ 
&\genfrac{}{}{0pt}{1}{\text{apply}}{\rR\pr_*\rR\cHom(\,\cdot\,,\pl^*\beF(-3))} 
&\genfrac{}{}{0pt}{1}{\text{apply}}{\rR\pr_*\rR\cHom(\pl^*\beF^\vee(1),\,\cdot\,)}
&
}$}
}\label{eq:diag}
\end{m-eqn}
Hence the various homomorphism in~\eqref{eq:VW} are natural, functorial. 
\end{m-proposition}
The unusual displays of the top and rightmost exact sequences are such that the arrows induced in the next diagram are in normal position.

\begin{m-proof}
By applying the indicated functors, we obtain: 
\begin{m-eqn}{
\kern-3ex\scalebox{.8}{$
\xymatrix@R=2.25em@C=1.15em{
\mbox{$\underbrace{H^2(\beF(-3))}_{=\,V_\beF}\otimes\eO_\ppp(-1)$}
\ar@{^(->}[r]^-{\veps_\beF}\ar@{.}[d]|-{?}
&
\mbox{$\underbrace{
\cExt^2_{\pr}\big(\eI_\Delta,\pl^*\beF(-3)\big)\otimes\eO_\ppp(-1)}_{\kern-7ex=:\;\euf A}
$}
\ar@{->>}[r]\ar@{^(->}[d]^-{a_\beF^\vee}
&
\cExt^3_\pr(\eO_\Delta,\beF(-4))\!=\!\beF\!=\!\pr_*\pl^*(\beF_\Delta)
\ar@{^(->}[d]
\\ 
R^1\pr_*(\pl^*\beF(-1)\,{\otimes}\,\eS\,{\otimes}\,\pr^*\eO_\ppp(1) )
\ar@<-.5ex>@/_1ex/@{^(->}[r]
&
\ouset
{C^{\euf B}_{\beF}:= H^1(\ppp,\,\beF(-1)\otimes\Omega_\ppp^1(1)\,)\,\otimes\,\eO_\ppp}
{C^{\euf A}_{\beF}:= \Ext^2_\ppp(\Omega_\ppp^1(1),\beF(-3))\,\otimes\,\eO_\ppp}
{{\fbox{$\cong$}}}
\ar@/_1ex/@<-.5ex>@{->>}[r]^-{b_\beF}\ar@{->>}[d]
&
\mbox{$\overbrace{
R^1\pr_*\big(\eI_\Delta\otimes\pl^*\beF(-1)\big)\otimes 
\eO_\ppp(1)}^{\kern7ex=:\;\euf B}
$}
\ar@{->>}[d]^-{q_\beF}
\\ 
&
\cExt^2_\pr(\eS\otimes\pr^*\eO_\ppp(1),\pl^*\beF(-3) )
\ar@{.}[r]|-{?}
&
\mbox{$\underbrace{H^1(\beF(-1))}_{=\,W_\beF}\otimes\eO_\ppp(1)$}
}
$}
}\label{eq:beil}
\end{m-eqn}\vspace{-.5ex}
where $\cExt_\pr$ stands for the relative $\Ext$-functor. (For the middle row/column, one uses~\eqref{eq:eS} to deduce the cohomology vanishings involving $\eS$.)
The rightmost and top extensions are given by the identity elements of $\End(W_\beF)$ and $\End(V_\beF)$: they are determined by the (unique) extension $\eI_\Delta\hra\eO_{\ppp\times\ppp}\surj\eO_\Delta$ tensored by $\beF(-1)$ etc. Thus $\euf A, \euf B$ are isomorphic to $\eK_\beF$ and $\eQ_\beF$, respectively. The fact that $a,b$ are evaluation morphisms follow from the identities: 
\begin{longtable}{rl}
$H^0(\euf B)
=H^0(\ppp,\,R^1\pr_*(\eI_\Delta\otimes\pl^*\beF(-1))\otimes\eO(1))$
&
$=H^1(\ppp\times\ppp,\,\eI_\Delta\otimes(\beF(-1)\boxtimes\eO(1))\,)$
\\[.5ex] 
$=H^1(\ppp,\,\beF(-1)\otimes\pl_*(\eI_\Delta\otimes\pr^*\eO(1))\,)$
&
$=H^1(\ppp,\,\beF(-1)\otimes\Omega_\ppp^1(1))=C_\beF^{\euf B};$
\\[.5ex]
$H^0(\euf A^\vee)^\vee=\ldots
=\Ext^2_\ppp(\pl_*(\eI_\Delta\otimes\pr^*\eO(1)),\beF(-3))$
&
$=\Ext^2_\ppp(\Omega_\ppp^1(1),\beF(-3))=C_\beF^{\euf A}.$
\end{longtable}\vspace{-.5ex}
Now we check, respectively, the isomorphism of the entries in the leftmost column and the bottom row. The Euler sequence and~\eqref{eq:eS} yield: 
\begin{m-eqn}{
\scalebox{.9}{$
\begin{array}{l}
R^1\pr_*(\pl^*\beF(-1)\otimes\eS\otimes\pr^*\eO_\ppp(1) )\cong H^1(\beF{\otimes}\eT_\ppp(-3))\otimes\eO_\ppp(-1)\cong H^2(\beF(-3))\otimes\eO_\ppp(-1),
\\[1.5ex] 
\cExt^2_\pr(\eS\otimes\pr^*\eO_\ppp(1),\pl^*\beF(-3))\cong H^2(\Omega^1_\ppp(1)\otimes\beF(-2))\otimes\eO_\ppp(1)\cong H^1(\beF(-1))\otimes\eO_\ppp(1).
\end{array}
$}
}\label{eq:R12}
\end{m-eqn}
Therefore~\eqref{eq:beil} agrees with~\eqref{eq:VW}, up to isomorphism; the latter is determined by extending to the left the rightmost column, by using $\Ext^1(W_\beF(1),\eK_\beF)\srel{\cong}{\to}\Ext^1(W_\beF(1),\beF).$ 
\end{m-proof}

We need to explicitly determine the maps in~\eqref{eq:beil} and also the mysterious isomorphism between $C_\beF^{\euf A}, C_\beF^{\euf B}$ (this is necessarily so, by general considerations). 
Moreover, the left- and lowermost terms in~\eqref{eq:beil} are not equal --they are isomorphic~\eqref{eq:R12}-- which is confusing for doing cohomological computations. Things get straightened out by restricting~\eqref{eq:diag} to a general hyperplane $\cH\subset\ppp$. The verification of the following claim is straightforward.

\begin{m-lemma}\label{lm:hom2}
Let $\beF_\cH$ be an instanton on $\cH$. (Thus, in here, $\beF_\cH$ is not necessarily the restriction of some $\beF$ to $\cH$.) By applying the indicated functors to the diagram, 
\vskip0.5em\begin{m-eqn}{
\scalebox{.8}{$
\xymatrix@R=1em@C=1.5em{
\pr^*\eO_\cH(1)\ar@/^4ex/@{->>}[rr] &\,\eI_{\Delta_\cH}\otimes\pr^*\eO_\cH(1)\ar@{_(->}[l]&\eO_{\Delta_\cH}(1) 
&\genfrac{}{}{0pt}{1}{\text{apply}}{\rR\pr_*\rR\cHom(\,\cdot\,,\pl^*\beF_\cH(-2))} 
\\ 
\eO_\cH(-1)\boxtimes\eO_\cH(-1)\,\ar@{^(->}[r]
&\pl^*\Omega^1_\cH(1)\ar@{->>}[r]\ar@{->>}[u]
&\eI_{\Delta_\cH}\otimes\pr^*\eO_\cH(1)\ar@{^(->}[d] 
&\genfrac{}{}{0pt}{1}{\text{apply}}{\rR\pr_*\rR\cHom(\pl^*\beF_\cH^\vee(1),\,\cdot\,)}
\\ 
&\eO_\cH(-1)\boxtimes\eO_\cH(-1)\ar@{^(->}[u]
&\pr^*\eO_\cH(1)\ar@/_10ex/@{->>}[uu] 
& 
\\ 
&\genfrac{}{}{0pt}{1}{\text{apply}}{\rR\pr_*\rR\cHom(\,\cdot\,,\pl^*\beF_\cH(-2))} 
&\genfrac{}{}{0pt}{1}{\text{apply}}{\rR\pr_*\rR\cHom(\pl^*\beF_\cH^\vee(1),\,\cdot\,)}
&
}$}
}\label{eq:diag2}
\end{m-eqn}
\nit one obtains the Barth-Hulek display of the monad corresponding to $\beF_\cH$:
\begin{m-eqn}{
\kern-1ex\scalebox{.8}{$
\xymatrix@R=2.25em@C=1.5em{
\scalebox{.95}{$
\underbrace{H^1(\beF_\cH(-2))}_{=\,V_{\beF_\cH}}\otimes\eO_\cH(-1)
$}\ar@{^(->}[r]^-{\veps_{\beF_\cH}}\ar@{=}[d]
&
\mbox{$\underbrace{
\scalebox{.9}{$\cExt^1_{\pr}\big(\eI_{\Delta_\cH},\pl^*\beF(-2)\big)\otimes\eO_\cH(-1)$}}_{\kern-7ex=:\;\euf A_\cH}
$}
\ar@{->>}[r]\ar@{^(->}[d]^-{a_{\beF_\cH}^\vee}
&
\scalebox{.9}{$\cExt^2_\pr(\eO_{\Delta_\cH},\beF_\cH(-3))\!=\!\beF_\cH\!=\!\pr_*\pl^*(\beF_{\Delta_\cH})$}
\ar@{^(->}[d]
\\ 
\scalebox{.9}{$H^1(\beF_\cH(-2))\otimes\eO_\cH(-1)$}
\ar@<-.5ex>@/_1ex/@{^(->}[r]
&
\ouset{C^{\euf B}_{\beF_\cH}:= H^1(\beF_\cH(-1)\otimes\Omega_\cH^1(1)\,)\,\otimes\,\eO_\cH}
{C^{\euf A}_{\beF_\cH}:= H^1(\beF_\cH(-2)\otimes\eT_\cH^1(-1))\,\otimes\,\eO_\cH}
{{\fbox{$=$}}}
\ar@/_1ex/@<-.5ex>@{->>}[r]^-{b_{\beF_\cH}}\ar@{->>}[d]
&
\mbox{$\overbrace{
\scalebox{.9}{$R^1\pr_*\big(\eI_{\Delta_\cH}\otimes\pl^*\beF_\cH(-1)\big) \otimes \eO_\cH(1)$}}^{\kern7ex=:\;\euf B_\cH}
$}
\ar@{->>}[d]^-{q_{\beF_\cH}}
\\ 
&
\scalebox{.9}{$H^1(\beF_\cH(-1))\otimes\eO_\cH(1)$}\ar@{=}[r]
&
\scalebox{.9}{$
\underbrace{H^1(\beF_\cH(-1))}_{=\,W_{\beF_\cH}}\otimes\eO_\cH(1)
$}
}
$}
}\label{eq:beil2}
\end{m-eqn}
\end{m-lemma}
We remark that the dimensional reduction transforms the dotted isomorphisms in~\eqref{eq:beil} into equalities. For a $2$-plane $\cH\subset\ppp$, the resolution of the diagonal $\Delta_\cH\subset \cH\times \cH$ is 
$$
\xymatrix@C=2.5em{
\ouset{}{2}{{\mbox{$\bigwedge$}}}
\bigl(\Omega_\cH^1(1)\boxtimes\eO_\cH(-1)\bigr)
=
\underbrace{\Omega_\cH^2(2)\boxtimes\eO_\cH(-2)}_{=\eO_\cH(-1)\boxtimes\eO_\cH(-2)}
\ar@{^(->}[r]^-{s_\cH\ort}&
\Omega_\cH^1(1)\boxtimes\eO_\cH(-1)
\ar@{->>}[r]^-{s_\cH\ort}&
\eI_{\Delta_\cH},
}
$$
where the homomorphisms are contractions with the identity  $s_\cH\in\End(H^0(\eO_\cH(1)))$, the restriction of $s$ to $\cH$. Note also that the tangent sequence splits: 
$$\,
{\eT_\ppp(-1)}\rst_\cH=\eT_\cH(-1)\oplus\eO_\cH,\;\;
{\Omega_\ppp^1(1)}\rst_\cH=\Omega_\cH^1(1)\oplus\eO_\cH.
$$ 

\begin{m-lemma}
The isomorphism $C_\beF^{\euf A}\to C_\beF^{\euf B}$ fits into the following commutative diagram: 
\begin{align*}
\scalebox{.8}{\kern-2ex$
\xymatrix@R=3em@C=1.5em{
V_\beF\!\otimes\!\eO_\ppp(-1)
\ar[d]^-{s\otimes}
\ar[r]^-{\del^{-1}\otimes\res_\cH}_-{\cong\otimes\res_\cH}
\ar@/_9ex/[dd]|(.25){a_\beF^\vee\veps_\beF}
&
H^1(\beF_\cH(-2))\!\otimes\!\eO_\cH(-1)
\ar[d]^-{s_\cH\otimes}
\ar@/_6ex/[dd]|(.25){a_{\beF_\cH}^\vee\veps_{\beF_\cH}}
\ar[r]^{0}_{\rm zero\,map!}
&
H^1(\beF_\cH(-1))\!\otimes\!\eO_\cH(1)
&
W_\beF\!\otimes\!\eO_\ppp(1)\ar[l]_-{\res_\cH\otimes\res_\cH}^-{\cong\otimes\res_\cH}
\\ 
\scalebox{.9}{$V_\beF\!\otimes\! H^0(\eT_\ppp(-1))$} 
\ar[d]^-\ev 
&
\scalebox{.9}{$
V_{\beF_\cH}H^0(\eT_\cH(-1))\ar[d]^\ev H^0(\eO_\cH(1))\eO_\cH(-1)$}
\ar[d]
&
\scalebox{.9}{$
W_{\beF_\cH}\!\otimes\! H^0(\eO_\cH(1))$}
\ar[u]
&
\scalebox{.9}{$W_\beF\!\otimes\! H^0(\eO_\ppp(1))$}
\ar[u] 
\\ 
\scalebox{.9}{$
\underbrace{H^2(\beF(-3)\otimes\eT_\ppp(-1))}_{=\,C_\beF^{\euf A}}
$}
\ar[r]^-{\del^{-1}\otimes\res_\cH}_-\cong
\ar@/_7ex/[rrr]|-{\;\cong\;}
&
\scalebox{.9}{$
\underbrace{H^1(\beF_\cH(-2)\otimes\eT_\cH(-1))}_{=\,C_{\beF_\cH}^{\euf A}}
$}\ar@{=}[r]
&
\scalebox{.9}{$
\underbrace{H^1(\beF_\cH(-1)\otimes\Omega_\cH^1(1))}_{=\,C_{\beF_\cH}^{\euf B}}
$}
\ar[u]^-{s_\cH\ort}
\ar@/_10ex/[uu]|(.3){q_{\beF_\cH} b_{\beF_\cH}}
&
\scalebox{.9}{$
\underbrace{H^1(\beF(-1)\otimes\Omega_\ppp^1(1))}_{=\,C_\beF^{\euf B}}
$} 
\ar[u]^-{s\ort}
\ar[l]_-{\res_\cH}^-\cong\ar@/_10ex/[uu]|(.3){q_\beF b_\beF}
}$}
\tag{CAB}\label{eq:cab}
\end{align*}
\end{m-lemma}
(At the $(2,2)$-entry we suppressed the `$\otimes$', due to lack of space.)
\begin{m-proof}
The curved arrows are obtained by applying $\cHom$-functors to $\pl^*\Omega_\ppp^1\to\pr^*\eO_\ppp(1)$ and to its restriction to $\cH$. We prove the factorization for $a_{\beF}^\vee\veps_{\beF}$, the other cases are similar. Note that $a_{\beF}^\vee\veps_{\beF}$ is obtained by applying $R^1\pr_*(s)$ to the homomorphism, dual to Beilinson's map,  $\eO_\ppp\boxtimes\eO_\ppp(-1)\srel{s\,\otimes}{\lar}\eT_\ppp(-1)\boxtimes\eO_\ppp$, tensored by $\pl^*\beF(-2)$. By pairing first with the $H^0(\eO_\ppp(1))$-component, we obtain the factorization. 

Now we turn our attention to the horizontal arrows. Let $\cH\subset\ppp$ be a general hyperplane. 
\begin{flushleft}
\begin{minipage}[t]{0.5\textwidth}
We tensor by $\beF(-1)$ the commutative diagram on the right and deduce that the restriction $$H^1(\beF\otimes\Omega^1_\ppp)\to H^1(\beF_\cH\otimes\Omega^1_\cH)$$ is an isomorphism.
\end{minipage}\quad 
\begin{scalebox}{.9}{
$\xymatrix@R=1.5em@C=1.5em{
\Ker\;\ar@{^(->}[r]\ar@{^(->}[d]
&\eO_\ppp\oplus\eO_\ppp(-1)^{\oplus 3}\ar@{^(->}[d]\ar@{->>}[r]
&\eO_\ppp\ar@{^(->}[d]
\\ 
\Omega^1_\ppp(1)\ar@{^(->}[r]\ar@{->>}[d]
&\eO_\ppp\oplus\eO_\ppp^{\oplus 3}\ar@{->>}[d]\ar@{->>}[r]
&\eO_\ppp(1)\ar@{->>}[d]
\\ 
\Omega^1_\cH(1)\ar@{^(->}[r]&\eO_\cH^{\oplus 3}\ar@{->>}[r]&\eO_\cH(1)
}$
}\end{scalebox}
\end{flushleft}
A similar argument shows that $H^1(\beF_\cH\otimes\eT_\cH(-3))\srel{\del}{\to} H^2(\beF\otimes\eT_\ppp(-4))$ is an isomorphism too. (Or simply take the Serre-dual and reduce to the previous case.) 

But $\eT_\cH(-1)=\Omega^1_\cH(2)$, so $H^1(\beF(-1)\otimes\Omega^1_\ppp(1))$, $H^2(\beF(-3)\otimes\eT_\ppp(-1))$ are both isomorphic to $H^1(\beF_\cH\otimes\Omega^1_\cH)$, when restricted to $\cH$.
\end{m-proof}

One might wonder what's the use of this diagram, since the composed `down-then-up' homomorphism vanishes. (We go from the first to the last entry of the exact sequence $V_\beF(-1)\to\dots W_\beF(1)$.)  We'll see that the $\eO(\pm1)$-terms in the middle row, causing the vanishing of the evaluation maps, are absorbed into the cohomology of another instanton which enters into the picture. Thus, we'll definitely deal with non-zero maps.


\section{Determining homomorphisms}\label{sct:hom} 

In the sequel, we keep in mind that cohomology classes are represented by {\v C}ech cocycles, which are genuine sections over open subsets. This is useful for understanding the effect of various homomorphisms. Cocycles are commonly denoted by $\cal Z^\bullet(\,\cdot\,)$. 

Let $\beF, \beG$ be instanton vector bundles, of \emph{possibly different ranks and charges}! We consider the display~\eqref{eq:beil} for $\beF$, and let $\delta_{\beF}$ be the boundary map in cohomology, corresponding to the top line. The tensor product with $\beG(-2)$ yields the diagram: 
\begin{align*}
\kern-3ex\scalebox{.9}{$
\xymatrix@R=3em@C=1em{
&&
\null\hspace{4em} 
W_{\beF}\otimes H^1(\beG(-1))
\ar[d]^-{\delta_{\beF}\otimes\bone_\beG}_-\cong
\ar[dl]_-{\kz_{\beF}\otimes\kz_\beG}^(.37)\cong 
&
\\ 
0\to  H^1(\beF\otimes\beG(-2))\ar[r] 
& 
V_{\beF}\otimes H^2(\beG(-3)) 
\ar[r]_-{ H^2(\veps_{\beF}\otimes\bone_\beG)}&
H^2(\eK_{\beF}\otimes\beG(-2))\ar[r] 
& H^2((\beF\otimes\beG)(-2))\to0
}
$}\tag{KZ}\label{eq:kz}
\end{align*}

\begin{m-proposition}\label{prop:=} 
{\rm(i)} The triangle in the diagram above is commutative that is, 
\begin{m-eqn}{
(\delta_{\beF}\otimes\bone_\beG)^{-1}\circ H^2(\veps_{\beF}\otimes\bone_\beG)
=\kz_\beF^{-1}\otimes\kz_\beG^{-1}.
}\label{eq:comm}
\end{m-eqn}
\nit{\rm(ii)} The tensor product of two mathematical instantons is still a mathematical instanton:
$$
H^1(\beF\otimes\beG(-2))=H^2(\beF\otimes\beG(-2))=0.
$$
\end{m-proposition} 

\begin{m-remark}\label{rk:tpym}
The Atiyah-Hitchin-Singer correspondence~\cite[Theorem 5.2]{ahs} implies that the tensor product of two Yang-Mills instantons on $\mbb R^4\cup\{\infty\}$ is still a Yang-Mills instanton, so the $H^1$- and $H^2$-cohomologies of the $(-2)$-twist vanish. Our result generalizes this property and the proof is algebraic. We discuss this topic and the relevance for physics in Section~\ref{sct:roots}. 
\end{m-remark}

The proposition has categorical interpretation (see~\cite{egno} for the definitions). Let $\IVB$ be the category whose objects are instantons on $\ppp$, of variable rank and charge (they are automatically semi-stable); let $\IVB^{\rm ps}$ be the subcategory formed by poly-stable bundles, finite direct sums of stable objects. 

\begin{m-theorem}\label{thm:categ}
{\rm(i)} The direct sum and tensor product define morphisms between the moduli spaces of mathematical instantons:
$$
\begin{array}{l}
\oplus:\mi_\ppp(r';n')\times\mi_\ppp(r'';n'')\to\mi_\ppp(r'+r'';n'+n''),\\ 
\otimes:\mi_\ppp(r';n')\times\mi_\ppp(r'';n'')\to\mi_\ppp(r'r'';r''n'+r'n'')
\end{array}
$$
Thus Schur powers (representations) preserve mathematical instantons. 

\nit{\rm(ii)} 
$(\IVB,\otimes)$ is a symmetric monoidal category, with unit $\eO_\ppp$. The opposite category $\IVB^{\rm op}$ is equivalent to $\IVB$ by the duality functor $\IVB^{\rm op}\to\IVB,\;\beF\to\beF^\vee$.  
 
\nit{\rm(iii)} $(\IVB^{\rm ps},\otimes)$ is actually a multi-tensor subcategory. 
\end{m-theorem}
In the last statement, one must consider poly-stable objects, as tensor products of stable vector bundles may be decomposable: e.g. $\beF\otimes\beF=Sym^2(\beF)\oplus\overset{2}{\bigwedge}\beF.$ One may view $\IVB^{\rm ps}$ as the quotient of $\IVB$ by the equivalence relation determined by the Jordan-H\"older filtration. 

\begin{m-proof}
The statements follow from the previous proposition and the fact that the tensor product of two semi-stable (resp. poly-stable) vector bundles is still semi-stable (resp. poly-stable). Over $\mbb C$, this is a consequence of the Kobayashi-Hitchin correspondence.
\end{m-proof}


\subsection{Proving~{\rm(KZ)}}
One may rephrase the Proposition as follows: 
$$
(\delta_{\beF}\otimes\bone_\beG)^{-1}\circ H^2(\veps_{\beF}\otimes\bone_\beG):V_{\beF}\otimes V_\beG\to W_{\beF}\otimes W_\beG
$$ 
is a natural assignment from instanton bundles to homomorphisms between vector spaces, so it's natural to ask what is this map. 
On both sides, $\beF,\beG$ are independent, and this leads to the idea of analysing the effect on $\beF$ and $\beG$ separately. It's quite confusing that, although the final (desired) conclusion is completely symmetric in the entries, this is obtained by composing `very asymmetric' terms. Note two re-arrangements, which justify the appearance of restrictions to $2$-planes in $\ppp$ in the sequel: 
\begin{m-eqn}{
V_\beF\otimes\beG(-3)=V_\beF(-1)\otimes\beG(-2)\quad\text{and}\quad W_\beF(1)\otimes\beG(-2)=W_\beF\otimes\beG(-1).
}\label{eq:*/}
\end{m-eqn}
They are necessary to apply $\veps_\beF,\delta_\beF$, and correspond to division (for $V_\beF,W_\beF$) by a linear equation and multiplication (for $\beG$) by the same linear factor, respectively. Our reasoning involves three steps: first, we analyse the effect of the homomorphism $\veps_\beF\otimes\bone_\eG$; second, we analyse $\delta_\beF\otimes\bone_\eG$; finally, we compose the two maps. 

\begin{m-lemma}\label{lm:dta}
The homomorphism $\phi$ in~\eqref{eq:dta} below has the form   
$\phi=\chi_\phi\otimes\bone_{H^2(\beG(-3))},$ where the $\beF$-component is $\chi_\phi\in\Hom(V_\beF,W_\beF).$
\end{m-lemma}

\begin{m-proof} 
Let $\lda$ be a trivialising line for $\beG$ and $\cD, \cH$ be two planes containing it; the restrictions $\beG_\cD,\beG_\cH$ are automatically semi-stable. We claim that the following diagram is commutative: 
\begin{m-eqn}{
\scalebox{.85}{$
\xymatrix@R=3em@C=1.25em{
V_{\beF_\cH}\otimes H^1(\beG(-2)_\cH)
\ar@/_14ex/[ddd]|(.5)\cong^(.57){\del\otimes\del}
\ar@{^(->}[r]^-{\veps_{\beF_\cH}}
\ar@/_5.5ex/@{-->}[rrr]^(.5){\;\phi_\cH\;}
\ar@/^5ex/[rr]|(.5){\;s_\cH\ort\;}
&
H^1(\eK_{\beF_\cH}\otimes\beG(-1)_\cH)\ar@{^(->}[r]
&
C_{\beF_\cH}^{\euf A}\otimes H^1(\beG(-1)_\cH)\ar@{=}[r]^{\rm incl.}_{\lda\subset \cH}
&
C_{\beF_\cH}^{\euf A}\otimes H^1(\beG(-2)_\cH)\ar@{->>}[dd]|-\cong^(.57){\del\otimes\del}
\\ 
V_\beF\otimes H^1(\beG(-2)_{\cD\cup \cH})\ar@{->>}[u]\ar[d]^-\cong\ar[r]
&
H^1(\eK_\beF\otimes\beG(-1)_{\cD\cup \cH})\ar[d]\ar[r]\ar[u]
&
C_\beF^\eA\otimes H^1(\beG(-1)_{\cD\cup \cH})\ar@{->>}[u]\ar@{->>}[dr]
&
\\ 
V_\beF\otimes H^2(\beG(-4))\ar[d]^-{\gamma_\beG}\ar[r]^-{\veps_\beF}
&
H^2(\eK_\beF(1)\otimes\beG(-4))\ar@{->>}[d]_(.35){\mu_\cH}\ar[rr]|(.4){\;\cong\;}^(.4){a_\beF^\vee}
&&
C_\beF^\eA\otimes H^2(\beG(-3))
\\ 
V_\beF\otimes H^2(\beG(-3))\ar[r]_-{H^2(\veps_\beF\otimes\bone_\beG)}\ar@/^1.5ex/@{-->}[rrru]_(.5){\;\phi\;}
&
H^2(\eK_\beF(1)\otimes\beG(-3))
&
&
}
$}
}\label{eq:dta}
\end{m-eqn}
Indeed, except the leftmost column --it's actually a square-- and the top-right square, all the arrows are natural homomorphisms in the display of the monad (horizontally), and restrictions (vertically). The leftmost column is obtained by twisting the third diagram~\eqref{eq:IDH} with $\eG(-2)$. The top-right square involving the inclusion $\lda\subset \cH$ is obtained by tensoring with $\beG(-1)$ the commutative diagrams: 
$$
\begin{array}{ll|ll}
\xymatrix@R=1.5em@C=1.5em{\eI_{\cD\cup \cH}\ar@{=}[d]\ar[r]&\eI_\cD\ar[d]\ar[r]&\eO_\cH(-\lda)\ar[d]
\\ 
\eI_{\cD\cup \cH}\ar[r]&\eO_\ppp\ar[r]&\eO_{\cD\cup \cH},}
&&&
\xymatrix@R=1.5em@C=1.5em{
\eO_\cH(-\lda)\ar[r]\ar[d]&\eO_\cH\ar@{=}[d] \\ 
\eO_{\cD\cup\cH}\ar[r]&\eO_\cH.}
\end{array}
$$
We declare that $H^1(\beG_\cH(-2))$ and $H^1(\beG_\cH(-1))$ are equal, because the isomorphism between them is determined by the inclusion $\eO_\cH(-1)\subset\eO_\cH$. The second row `interpolates' between the first and the third. It allows defining the dashed homomorphism $\phi$ by following the left-top-right path, involving restrictions to $\cH$. 

Let us prove that $\phi$ acts as the identity on $H^2(\beG(-3))$; equivalently, the $\beG$-component of $\phi_\cH$ is the identity of $H^1(\beG_\cH(-2))$. 
The composition of the first two top arrows --see~\eqref{eq:cab}-- is the pairing with $s_\cH\in\Hom(V_{\beF_\cH},C^\eA_{\beF_\cH})\otimes H^0(\eO_\cH(1))$, and it factorizes: 
$$
V_{\beF_\cH}\otimes H^1(\beG_\cH(-2))\srel{s_\cH}{\lar} C_{\beF_\cH}^\eA\otimes H^0(\eO_\cH(1))\otimes H^1(\beG_\cH(-2))
\srel{\bone_C\otimes\ev}{\lar} C_{\beF_\cH}^\eA\otimes H^1(\beG_\cH(-1)).
$$ 
In coordinates $\zeta_{\cH,j},\,j=0,\dots,2$ on $\cH$, we have 
$s_\cH=\ouset{j=0}{2}{\sum}c_{\cH,j}\zeta_{\cH,j},\; c_{\cH,j}\in\Hom(V_{\beF_\cH},C^\eA_{\beF_\cH}).$ We use $\zeta_{\cH,j}$ for defining inclusions $\beG_\cH(-2)\to\beG_\cH(-1)$; they induce equality in the degree-$1$ cohomology. Since we used the same linear forms both in the evaluation and for the inclusion maps, their composition is the identity of $H^1(\beG_\cH(-2))$. We note that the $\beF_\cH$-component of $\phi_\cH$ is $\chi_\phi=\ouset{j=0}{2}{\sum}c_{\cH,j}$, and it's independent of $\beG$.

It may be illuminating to give a second proof when, for general $\cH$, the restriction $\beG_\cH$ is stable ($\beG$ is stable, too), thus $\gamma_\beG$ is surjective.  We directly apply the division-multiplication trick~\eqref{eq:*/}: a lifting of an element in $H^2(\beG(-3))$ to $H^2(\beG(-4))$ amounts to dividing the corresponding cocycle by a linear equation. (The surjectivity of the arrow ensures that such division makes sense, and the commutativity of the diagram implies that the result is independent of the lifting.) The homomorphism $\veps_\beF$ has the form $\ouset{j=0}{3}{\sum} c_j\zeta_j$, where $\zeta_0,\dots,\zeta_3$ are coordinates on $\ppp$ and $c_j\in\Hom( V_\beF, C_\beF)$. By following the third row, we see that $\phi$ acts on $v\otimes h\in V_\beF\otimes H^2(\beG(-3))$ as follows: for $j=0,\dots,3$, there is a representative $\tld h_j\in\cal Z^2(\beG(-3))$ of $h$, such that the quotient $\frac{\tld h_j}{\zeta_j}\in\cal Z^2(\beG(-4))$ is well-defined, so we have 
$$
\phi(v\otimes h)=\text{the cohomology class defined by}\;\bigg[\ouset{j=0}{3}{\sum} c_j(v)\zeta_j\otimes\frac{\tld h_j}{\zeta_j}\bigg]
=\biggl[\ouset{j=0}{3}{\sum} c_j(v)\biggr]\otimes h.
$$
Thus, the linear factor required for lifting $h$ to $H^2(\beG(-4))$ cancels out by applying $\veps_\beF$. 
\end{m-proof}

\begin{m-lemma}\label{lm:gma}
The homomorphism $\psi$ in~\eqref{eq:gma} below is of the form  
$\psi=\chi_\psi\otimes\kz_\beG^{-1},$ where the $\beF$-component is $\chi_\psi\in\Hom(C^\eA_\beF,W_\beF).$
\end{m-lemma}

\begin{m-proof} The tensor product of the sequences 
$$\eK_\beF(-1)\hra \bigl(C_\beF^\eA(-1)\srel{\cong}{\lar} C_\beF^\eB(-1)\bigr)\surj W_\beF\otimes\eO_\pp
\quad{\rm and}\quad
\beG(-2)\hra\beG(-1)\surj\beG(-1)_\cH$$ 
yield the diagram: 
$$
\xymatrix@R=1.5em@C=3em{
\eK_\beF\otimes\beG(-3)\ar@{^(->}[r]^-{\mu_\cH}\ar@{^(->}[d]
&
\eK_\beF\otimes\beG(-2)\ar@{^(->}[d]\ar@{->>}[r]
&
\eK_\beF\otimes\beG(-2)_\cH\ar@{^(->}[d]\ar@{-->}[d]
\\ 
C_\beF^\eA\otimes\beG(-3)\ar@{^(->}[r]\ar@{->>}[d]
&
C_\beF^\eA\otimes\beG(-2)\ar@{->>}[d]\ar@{->>}[r]\ar@{..>}[dr]|{\quad m\quad}
&
C_\beF^\eA\otimes\beG(-2)_\cH\ar@{->>}[d]
\\ 
W_\beF\otimes\beG(-2)\ar@{^(->}[r]
&
W_\beF\otimes\beG(-1)\ar@{->>}[r]
&
W_\beF\otimes\beG(-1)_\cH.
}
$$
Then $\euf M:=\Ker(m)$ satisfies $\,W_\beF\otimes H^1(\beG(-1)_\cH)\srel{\cong}{\to} H^2(\euf M)\,$, and it fits into: 
$$
\begin{array}{c|c|c}
\kern-1ex\scalebox{.62}{$\xymatrix@R=1.5em@C=1.25em{
\eK_\beF\otimes\beG(-3)\ar@{^(->}[r]\ar@{^(->}[d]
&
\eK_\beF\otimes\beG(-2)\ar@{^(->}[d]\ar@{->>}[r]
&
\eK_\beF\otimes\beG(-2)_\cH\ar@{=}[d]
\\ 
C_\beF^\eA\otimes\beG(-3)\ar@{^(->}[r]\ar@{->>}[d]
&
\euf M\ar@{->>}[r]\ar@{->>}[d]
&
\eK_\beF\otimes\beG(-2)_\cH
\\ 
W_\beF\otimes\beG(-2)\ar@{=}[r]
&
W_\beF\otimes\beG(-2)
&
}$}\kern-1ex
& 
\kern-1ex\scalebox{.62}{$\xymatrix@R=1.5em@C=1.25em{
&
&
\eK_\beF\otimes\beG(-2)_\cH\ar@{^(->}[d]
\\ 
C_\beF^\eA\otimes\beG(-3)\ar@{^(->}[r]\ar@{^(->}[d]
&
C_\beF^\eA\otimes\beG(-2)\ar@{=}[d]\ar@{->>}[r]
&
C_\beF^\eA\otimes\beG(-2)_\cH\ar@{->>}[d]
\\ 
\euf M\ar@{^(->}[r]\ar@{->>}[d]
&
C_\beF^\eA\otimes\beG(-2)\ar@{->>}[r]
&
W_\beF\otimes\beG(-1)_\cH
\\ 
\eK_\beF\otimes\beG(-2)_\cH
}$}\kern-1ex
& 
\kern-1ex\scalebox{.62}{$\xymatrix@R=1.5em@C=1.25em{
&
&
W_\beF\otimes\beG(-2)\ar@{^(->}[d]
\\ 
\eK_\beF\otimes\beG(-2)\ar@{^(->}[r]\ar@{^(->}[d]
&
C_\beF^\eA\otimes\beG(-2)\ar@{=}[d]\ar@{->>}[r]
&
W_\beF\otimes\beG(-1)\ar@{->>}[d]
\\ 
\euf M\ar@{^(->}[r]\ar@{->>}[d]
&
C_\beF^\eA\otimes\beG(-2)\ar@{->>}[r]
&
W_\beF\otimes\beG(-1)_\cH
\\ 
W_\beF\otimes\beG(-2)
}$}
\\ \rm (I)&\rm (II)&\rm (III)
\end{array}
$$
They imply the commutativity of the following diagram: 
\begin{m-eqn}{
\scalebox{.9}{$
\xymatrix@R=3em@C=3em{
H^2(\eK_\beF\otimes\beG(-3))\ar@{->>}[r]^-{\mu_\cH}\ar[d]^-\cong_-{a_\beF^\vee}\ar@{}[dr]|{\rm(I)}
&
H^2(\eK_\beF\otimes\beG(-2))\ar@{->>}[d]_(.35)\cong\ar@{}[dr]|{\rm(III)}
&
W_\beF\otimes H^1(\beG(-1))\ar[l]|(.5){\;\cong\;}_-{\delta_\beF\otimes\bone_\beG}\ar[d]^-\cong_-{\res_\cH}
\\ 
C_\beF^\eA\otimes H^2(\beG(-3))\ar@{->>}[r]\ar@/_0.25ex/@{-->}[urr]|(.37){\;\psi\;}\ar@{}[dr]|{\rm(II)}
&
H^2(\euf M)
&
W_{\beF_\cH}\otimes H^1(\beG(-1)_\cH)\ar[l]_-\cong
\\ 
C_{\beF_\cH}^\eA\otimes H^1(\beG(-2)_\cH)\ar@{->>}[u]_(.45)\cong^(.45)\del\ar[r]^-\cong
&
C_{\beF_\cH}^\eB\otimes H^1(\beG(-2)_\cH)\ar@{->>}[r]^-{s_\cH\ort}
&
W_{\beF_\cH}\otimes H^1(\beG(-1)_\cH)\ar@{->>}[ul]|-\cong\ar@{=}[u]
&
}
$}
}\label{eq:gma}
\end{m-eqn}
By moving along the lower edges of the diagram, we see that the homomorphism $\psi$ is a tensor product: its $\beF$-component is a composition of various homomorphisms between cohomology groups of $C_\beF,W_\beF$, so it depends only on $\beF$. 

Concerning the $\beG$-component, we claim that it acts on $H^2(\beG(-3))$ as the inverse of the Koszul map. Recall --see~\eqref{eq:kz-1}-- that the restriction to $\cH$ of $\kz^{-1}$ is simply the inclusion $H^1(\beG_\cH(-2))\srel{\cong}{\lar}H^1(\beG_\cH(-1))$. This is precisely the homomorphism obtained by following the lower side of the diagram, where we contract with the $H^0(\eO_\cH(1))$-component of $s_\cH$; see the third column of~\eqref{eq:cab}.
\end{m-proof}

\begin{proof}(of Proposition~\ref{prop:=}) 
(i) The composed homomorphism on the left-hand side of~\eqref{eq:comm} equals $\psi\circ\phi$, so is a tensor product of two linear maps. It's $\beG$-component is $\kz_\beG^{-1}$. 

To determine the $\beF$-component of $\chi_\psi\circ\chi_\phi$, we return once more to~\eqref{eq:cab}, and we analyse the restrictions to $\cH$ that is, the middle columns. Let us look at the $(2,2)$-entry in there: the homomorphism $\chi_\phi$ sending $V_{\beF_\cH}$ to $C^\eA_{\beF_\cH}$ in diagram~\eqref{eq:dta} is the evaluation (pairing) map applied to $H^1(\beF_\cH(-2))\otimes H^0(\eT_\cH(-1))$; as explained in the proof of Lemma~\ref{lm:dta}, the other half of the entry, namely $H^0(\eO_\cH(1))\otimes\eO_\cH(-1)$, is absorbed by the $\beG$-component during the division-multiplication process which yields the identity of $H^1(\beG_\cH(-2))$. 

So we reach the diagram~\eqref{eq:gma} and $\chi_\psi$. The bottom arrows start with $C_{\beF_\cH}^\eA\srel{=}{\to}C_{\beF_\cH}^\eB$, whose effect is a twist by $\eO_\cH(1)$ of the $\beF_\cH(-2)$-component, which yields $\beF_\cH(-1)$; the composition of the remaining arrows correspond to the (upward, $q_{\beF_\cH}b_{\beF_\cH}$) homomorphism in the third column of~\eqref{eq:cab}. Clearly, they act as the identity on $\beF_\cH(-1)$; once more, the $H^0(\eO_\cH(1))$-component of $s_\cH$ --the second bottom arrow in~\eqref{eq:gma} is the contraction with $s_\cH$-- is absorbed by the isomorphism $H^1(\beG_\cH(-2))\to H^1(\beG_\cH(-1))$.

Overall, $\chi_\psi\circ\chi_\phi$ twists $\beF_\cH(-2)$ by $\eO_\cH(1)$ and yields $H^2(\beF_\cH(-2))\to H^1(\beF_\cH(-1))$. Once more, \eqref{eq:kz-1} shows that this is nothing but the inverse of the Koszul homomorphism. 

\smallskip(ii) Since $\beF, \beG$ are instanton bundles, their Koszul maps are isomorphisms. Thus the extremities of the exact sequence~\eqref{eq:kz} vanish.
\end{proof}


\section{Mathematical instantons on the projective space}\label{sct:main}

In this section we prove the irreducibility and rationality properties of $\mi_\ppp(r;n)$ stated in the Introduction. 
Our approach consists in restricting instantons to either to a union (wedge) of planes (intersecting along a line) or to a smooth quadric. 

\subsection{Restriction maps}\label{ssct:res}

We consider the following geometric objects, which are general for the indicated properties: 
\begin{itemize}[leftmargin=3ex]
\item 
a line $\lda\subset\ppp$ and two $2$-planes $\cD, \cH$ intersecting along $\lda$.
\item 
another line $\lda'$ which intersects $\lda$ and $\cQ\cong\mbb P^1\times\mbb P^1$ a quadric containing $\lda\cup\lda'$. 
\end{itemize}
The moduli space $\mi_\ppp(r;n)$ has finitely many irreducible components. Since the choices above were general, for any component $M'\subset\mi_\ppp(r;n)$, the restrictions to $\lda, \lda'$ of the generic instanton bundle $\beF'$ in $M'$ are trivializable. Therefore the restrictions $\beF'_{\cD}, \beF'_{\cH}, \beF'_{\cQ}$ are all semi-stable. Thus we obtain the restriction maps (unframed and framed versions, dashed arrows stand for rational maps): 
\begin{m-eqn}{
\begin{array}{lll}
\Theta_{\cD\cH}:\mi_\ppp(r;n)\dashto \mm_{\cD\cup \cH}(r;n),
&&\beF\mapsto\beF_{\cD\cup \cH}:=\beF\otimes\eO_{\cD\cup \cH},
\\ 
\Theta_{\cD\cH,\lda}:\mi_\ppp(r;n)_\lda\rar \mm_{\cD\cup \cH}(r;n)_\lda, 
&&
\\[1ex] 
\Theta_{\cQ}:\mi_\ppp(r;n)\dashto \mm_{\cQ}(r;2n),
&&\beF\mapsto\beF_{\cQ}:=\beF\otimes\eO_{\cQ},
\\ 
\Theta_{\cQ,\lda\cup\lda'}:\mi_\ppp(r;n)_{\lda\cup\lda'}\rar \mm_{\cQ}(r;2n)_{\lda\cup\lda'},&&
\end{array}
}\label{eq:theta}
\end{m-eqn}
whose domains of definition meets all the irreducible components of $\mi_\ppp(r;n)$. (Note that, for $\cQ$, the charge is $2n$ because the quadric has degree two. The semi-stability property is with respect to the polarization $[\lda']+c[\lda], c>2r(r-1)n$.)

\begin{m-proposition}\label{prop:EgI}
\begin{enumerate}[leftmargin=5ex]
\item 
For any $\beF\in\mi_\ppp(r;n)$, the following properties hold: 
\begin{enumerate}
\item\label{itii:egen} 
$H^1((\cEnd\beF)(-2))=H^2((\cEnd\beF)(-2))=0;$
\item 
$H^2(\cEnd(\beF))=0$, so its deformations are unobstructed. 
\item[] The expected dimension of $\mi_\ppp(r;n)$ is $4rn-r^2+1$.
\end{enumerate} 
\item 
The differential of $\Theta_{\cD\cH}, \Theta_\cQ$ are isomorphisms everywhere, so they are \'etale maps.  
\item 
Each irreducible component of $\mi_\ppp(r;n)$ has the expected dimension and the locus corresponding to stable bundles is dense.
\end{enumerate}
\end{m-proposition}

\begin{m-proof}
(i) Just replace $\beG=\beF^\vee$ in Proposition~\ref{prop:=}. For the second property, let $\cH'$ be a general plane (for $\beF$), take the long exact sequences in cohomology determined by 
$$
\eO_\ppp(-2)\hra\eO_\ppp(-1)\surj\eO_{\cH'}(-1),
\quad\eO_\ppp(-1)\hra\eO_\ppp\surj\eO_{\cH'},
$$ 
twisted by $\cEnd(\beF)$, and use the semi-stabilty of $\beF_{\cH'}$. For a stable $\beF$, one computes that $h^1(\cEnd(\beF))=-\chi(\cEnd(\beF))+h^0(\cEnd(\beF))$ is indeed $4rn-r^2+1$. 

\nit(ii) The differentials of $\Theta_{\cD\cH}$ and $\Theta_\cQ$ at $\beF$ are, respectively, the homomorphisms 
$$
H^1(\cEnd(\beF))\to H^1(\cEnd(\beF_{\cD\cup \cH})),\quad H^1(\cEnd(\beF))\to H^1(\cEnd(\beF_{\cQ})).
$$ 
The property~\ref{itii:egen} shows that they are indeed isomorphisms. 

\nit(iii) Since $\Theta_{\cD\cH}$ (resp. $\Theta_\cQ$) is \'etale, its restriction to each component of $\mi_\ppp(r;n)$ is dominant. Stable vector bundles are dense in $\mm_{\cD\cup \cH}(r;n)$, resp. $\mm_{\cQ}(r;2n)$ (see Lemma~\ref{lm:stable}), and $\beF$ on $\ppp$ is stable as soon as its restriction to $\cD\cup \cH$, resp. $\cQ$, is so. Thus stable bundles are dense; at such a point, $\mi_\ppp(r;n)$ is smooth and has the expected dimension.
\end{m-proof}


\subsection{Irreducibility and rationality}\label{ssct:irred}
The following are our main results. 
\begin{m-theorem}\label{thm:birtl}
The restriction maps~\eqref{eq:theta} are birational. Actually,  $\Theta_{\cD\cH,\lda}, \Theta_{\cQ,\lda\cup\lda'}$ are open immersions.
\end{m-theorem}

\begin{m-proof} 
It is enough to prove the statement for the non-framed morphisms. 
Proposition~\ref{prop:EgI} implies that, restricted to each irreducible component $M'\subset\mi_\ppp(r;n)$, the map $\Theta_{\cD\cH}$ dominates $\mm_{\cD\cup \cH}(r;n)$ and  $\Theta_\cQ$ dominates $\mm_{\cQ}(r;2n)$. Suppose there are two irreducible components. Then there are two non-isomorphic stable instantons $\beF,\beG$ which are mapped to the same point. We proved in Proposition~\ref{prop:=} that $H^1(\cHom(\beF,\beG)(-2))=0$, so the isomorphism between the restrictions (either to $\cD\cup\cH$ or $\cQ$) lifts to an isomorphism over $\ppp$. 
\end{m-proof}

The result is in the same vein as~\cite{atiy-2+4,dons}: physical (Yang-Mills) instantons on $\mbb{CP}^3$ correspond to framed bundles on $\mbb{CP}^2$, resp. $\mbb{CP}^1\times\mbb{CP}^1$. We elaborate on this in Section~\ref{sct:roots}.

\begin{m-theorem}\label{thm:irred}
$\mi_\ppp(r;n)$ and $\mi_\ppp(r;n)_{\lda\cup\lda'}$ are irreducible and rational. 
\end{m-theorem}

\begin{m-proof}
We know that $\Theta_{\cD\cH}$, $\Theta_\cQ$ are birational. The statement follows from the irreducibility and rationality of $\mm_\cQ(r;2n)$ and $\mm_\pp(r;m)_\lda$, respectively (cf. Theorem~\ref{thm:rtl-frame}).
\end{m-proof}

The results obtained in Section~\ref{sct:hirz} (cf.~\eqref{eq:mihz} and {\S}\ref{ssct:applic}) yield a description of the general mathematical instanton on $\ppp$. 

\begin{m-corollary}\label{cor:B} 
The general mathematical instanton $\beF$ on $\mbb P^3$ is uniquely determined: 
\begin{enumerate}
\item[\underbar{\textbf{either}}] 
$\bullet\,$ by its restrictions $(\beF',\beF'')$ to $2$-planes $\cD,\cH\cong\pp$ intersecting along the line $\lda$;
\\ and\\ 
$\bullet\,$ by the gluing data $\beF'_\lda\cong\eO_\lda^{\oplus r}\cong\beF''_\lda$ (up to diagonal $\PGL(r)$-action). 
\item[] The general element of $\mm_{\pp}(r;n)$ is the kernel of a surjective homomorphism: 
\begin{m-eqn}{
\begin{array}{l} 
 {\eI_{p}^a(a)}^{\oplus r-\rho}\oplus{\eI_{p}^{a+1}(a+1)}^{\oplus \rho} \lar\ouset{j=1}{n}{\bigoplus}\;\eO_{l_j}(1),
 \quad 
a:=\lfloor n/r\rfloor,\, \rho:=n-ar, 
\\[2ex] 
l_1,\dots,l_n\subset\mbb P^2\;\text{are distinct lines passing through}\;p\in\mbb P^2. 
\end{array}
}\label{eq:mip}
\end{m-eqn}
\medskip\item[\underbar{\textbf{or}}] 
$\bullet\,$ by its restriction to a general quadric $\cQ\cong\mbb P^1_{\rm left}\times\mbb P^1_{\rm right}$. 
\item[] The general element of $\mm_{\cQ}(r;2n)$ is the kernel of a surjective homomorphism: 
\begin{m-eqn}{
\begin{array}{l}
{\eO_{\mbb P^1_{\rm left}}(a')}^{\oplus r-\rho'}\oplus{\eO_{\mbb P^1_{\rm left}}(a'+1)}^{\oplus\rho'}
 \lar\ouset{j=1}{2n}{\bigoplus}\;\eO_{\{x_j\}\times\mbb P^1_{\rm right}}(1),
\\[2ex] 
a':=\lfloor 2n/r\rfloor,\, \rho':=2n-a'r,\quad 
x_1,\dots,x_{2n}\in\mbb P^1_{\rm left}\;\text{are distinct points.} 
\end{array}
}\label{eq:miq}
\end{m-eqn}
\end{enumerate}
\end{m-corollary}


\section{Framed vector bundles on Hirzebruch surfaces}\label{sct:hirz}

The irreducibility and rationality of $\mi_\ppp(r;n)$ follow, once we know that the restriction map $\Theta_{\cD\cH}$ (resp. $\Theta_\cQ$) is birational, from the analogous statements for the moduli of vector bundles on $\pp$ (resp. $\mbb P^1\times\mbb P^1$). 
Note that $\pp$ is the blow-down of the $1^{\rm st}$ Hirzebruch surface; a quadric in $\ppp$ is isomorphic to $\mbb P^1\times\mbb P^1$, the $0^{\rm th}$ Hirzebruch surface. So we can utilize the techniques~\cite[\S3]{hlc-sheaves}, where the first named author studied vector bundles on Hirzebruch surfaces. In here, the difference is that we were led to framed instantons on $\ppp$, which involve different group actions, requiring changes. Since rationality is a sensitive issue, we tried to make the presentation short, yet (almost) self-contained.


\subsection{General properties} 

Let 
$Y_\ell:=\mbb P\bigl(\eO_{\mbb P^1}\oplus\eO_{\mbb P^1}(-\ell)\bigr)$ 
be the $\ell^{\rm th}$ Hirzebruch surface and $Y_\ell\srel{\pi}{\to}\mbb P^1$ the natural projection. 
We denote by $\eO_\pi(1)$ the relatively ample line bundle of $Y_\ell$,  $\Lambda:=\mbb P\bigl(\eO_{\mbb P^1}\oplus0\bigr)$ the $(-\ell)$-curve, and $l$ the general fibre of $\pi$; we have $[\eO_\pi(1)]=[\Lambda]+\ell[l]$. 

For integers $m\ges r\ges2$, we consider the polarization $L_c:=[\eO_\pi(1)]+c[l],\,c>mr(r-1),$ 
and the corresponding moduli space $\bar\mm_{Y_\ell}(r;m)$ of rank-$r$ torsion free sheaves $Y_\ell$, with $c_1=0$ and $c_2=m$. We denote $\mm_{Y_\ell}(r;m), \mm_{Y_\ell}(r;m)^{\svb}$ the open loci corresponding to vector bundles, resp. stable vector bundles.

\begin{m-lemma}\label{lm:stable}
For $m\ges r$, the generic vector bundle $\cV\in\mm_{Y_\ell}(r;m)$ is stable, the locus corresponding of stable bundles in dense.
\end{m-lemma}

\begin{m-proof}
Otherwise, the last term $\cV'$ of the Jordan-H\"older filtration of $\cV$ is a proper, saturated, semi-stable subsheaf, $\deg(\cV')=0$, it's reflexive, so locally free; $\cV'':=\cV/\cV'$ is torsion free, stable, $\deg(\cV'')=0$. 
Let $r':={\rm rank}(\cV')$, $m':=c_2(\cV')$, similarly for $\cV''$; Bogomolov's inequality yields $0\les m'\les m$. Note that $-h^1(\cV'') =\chi(\cV'')=r''-m''$, so $r''\les m''$. As $\cV$ is generic, its deformations are exhausted by deformations of $\cV', \cV''$ and extensions between them. 

We claim that, to the contrary, the following inequality holds true (The lower case $\ext,\dots,$ stand for the dimensions of the $\Ext,\dots$, respectively.):
 
\begin{longtable}{rl}
$1$&
$\les\ext^1(\cV,\cV)-[\ext^1(\cV',\cV')+\ext^1(\cV'',\cV'')-\ext^1(\cV'',\cV')]$
\\[1.5ex]
&
$=(2rn-r^2) - [(2r'n'-(r')^2)+(2r''n''-(r'')^2)+(r'n''+r''n'-r'r'')]$
\\[1.5ex]
&
$+[\dne(\cV)-(\dne(\cV')+\dne(\cV'')+\hom(\cV'',\cV'))]$
\\[1.5ex] &
$=Term_1+Term_2.$ 
\end{longtable}

A simple computation yields $Term_1=r'(n''-r'')+r''n'.$ Since $n''\ges r''$, it vanishes if and only if $n''=r'', n'=0$, which implies $n=n''=r''<r$, and this contradicts the hypothesis. 

We analyse the $Term_2$. The moduli space parametrizes \emph{classes} of sheaves up to Jordan-H{\"o}lder equivalence, so we can replace $\cV$ by  $JH(\cV)$ \textit{etc}; the situation becomes $\cV=\cV'\oplus\cV''$. The elements of $\End(\cV)$ have a block-form containing the other $\Hom$-spaces.
\end{m-proof}

\begin{thm-nono}{(cf.~\cite[2.8, 3.6]{hlc-sheaves})}
The following statements hold true:
\begin{enumerate}[leftmargin=5ex]
\item 
The $L_c$-semi-stability property (of torsion free sheaves) is independent of $c$ as above. If $\cV$ is $L_c$-semi-stable, its restriction to the general fibre $l$ is trivializable.  
\item 
The restrictions to both $\Lambda$ and $l$ of the generic $\cV\in\bar\mm_{Y_\ell}(r;m)$ are trivializable. 
\item[] (Note: the vector bundles involved in~\eqref{eq:theta} satisfy indeed this property.)
\item 
The generic $\cV$ is determined by an exact sequence of the form: 
\begin{m-eqn}{
\begin{array}{l}
0\to \pi^*\bbL\srel{\alpha}{\to}\cV\srel{\beta}{\to}
\cal S_\ux:=\ouset{j=1}{m}{\bigoplus}\eO_{\pi^{-1}(x_j)}(-1)
\to 0,\; \{x_1,\dots,x_m\}\subset\mbb P^1\;\textrm{distinct points,}
\\[3ex] 
\bbL:=\eO_{\mbb P^1}(-a)^{\oplus(r-\rho)}\oplus\eO_{\mbb P^1}(-a-1)^{\oplus\rho},\; 
a:=\lfloor m/r \rfloor,\;\rho:=m-ar.
\end{array}
}\label{eq:mihz}
\end{m-eqn}
\item 
$\bar\mm_{Y_\ell}(r;m)$ is irreducible, of dimension $2mr-r^2+1$. The exact sequences as above determine a unique maximal dimensional stratum.
\item 
The assignment $\cV\mt{\rm Supp}\,R^1\pi_*\cV(-\Lambda)$ defines a morphism $\bar\mm_{Y_\ell}(r;m)\srel{h}{\to}{\rm Hilb}^m_{\mbb P^1}\cong\mbb P^m$, whose generic fibre is $(2rm-r^2+1-m)$-dimensional. For $\cV$ as in~\eqref{eq:mihz}, we have $h(\cV)=\{x_1,\dots,x_m\}.$
\end{enumerate}
\end{thm-nono}

\begin{m-remark}\label{rk:bbr}
 The moduli space of framed sheaves on Hirzebruch surfaces, with framing along a section $l_\infty\in|\eO_\pi(1)|$ of $Y_\ell\srel{\pi}{\to}\mbb P^1$ was investigated by Bartocci-Bruzzo-Rava~\cite{bbr}. The authors allow arbitrary values for $r, m$, and possibly non-vanishing first Chern classes. 
 
 \begin{thm-nono}{(\cite[Theorem 3.4]{bbr})}
The moduli space of $l_\infty$-framed sheaves --for our purposes, we set $c_1=0, c_2=m$-- is smooth, irreducible, of dimension $2rm$, and it's fine (that is, it admits a universal Poincar\'e sheaf.) 

Any such framed sheaf is the cohomology of a complex of the form: 
 $$
 \eO_{Y_\ell}(-l)^{\oplus m}\to\eO_{Y_\ell}(l_\infty-l)^{\oplus m} \oplus\eO_{Y_\ell}^{\oplus(r+m)}\to\eO_{Y_\ell}(l_\infty)^{\oplus m}. 
 $$
 \end{thm-nono}
 
 Note that $l_\infty$ is a flat deformation of $\Lambda+\ell l$, so generic vector bundles $\cV$ as in Theorem(ii) above are trivializable along $l_\infty$, too. Since we are interested in birational properties, the result of Bartocci~\textit{et al.} yields the irreducibility of $\mm_{Y_\ell}(r;m)$ in our setup. 
 
The reason for working with the description~\eqref{eq:mihz} is that it is more economical compared to the detailed monad-type presentation, which involves the action of a large group. For successfully carrying out our computations, simplicity is essential. 

Lemma~\ref{lm:stable} shows that the condition $m\ges r$ ensures the density of the stable bundles, whose automorphism group consists of scalars, only. This is not longer true for $m<r$.  Indeed, for $\cV$ on $Y_\ell$, the Riemann-Roch formula yields: 
$$h^1(\cV)=(m-r)+h^0(\cV)\quad\Rightarrow\quad h^0(\cV)\ges r-m.$$ 
Thus every $\cV\in\mm_{Y_\ell}(r;m)$ admits non-trivial sections, so it's properly semi-stable.
\end{m-remark}


\subsection{The extension vector bundle} 
The explicit form of the general vector bundle determines a quotient description of $M_{Y_\ell}(r;m)^\svb$ and yields almost explicit coordinates on it. 

\subsubsection{The absolute case} 
Fix $0,\infty\in\mbb P^1$; $\mbb A^1=\mbb P^1\setminus\{\infty\}$ is the affine line. Let $\cA\subset\mbb A^m={\rm Hilb}^m(\mbb A^1)\subset\mbb P^m={\rm Hilb}^m(\mbb P^1)$ be the open locus of $m$-tuples $\ux=\{x_1,\dots,x_m\}$ consisting of distinct points on $\mbb A^1$. 
For $\ux\in\cA$, the extensions \textit{classes}~\eqref{eq:mihz} are parametrized by  
$$
E_{\ux}:=\ouset{i=1}{m}{\bigoplus}\;\mbb L_{x_i}
{\otimes}\,
\overbrace{\Gamma(\eO_{\pi^{-1}(x_i)}(1))}^{\cong\;\Bbbk^2}
=\Big(\mbb L\otimes\pi_*\eO_\pi(1)\Big)\otimes\eO_{\ux}, 
\quad\dim E_{\ux}=2mr.
$$
An element $e_\ux\in E_\ux$ determines~\eqref{eq:mihz}, equivalently, the dual form:
\begin{m-eqn}{
0\to\cV^\vee=\Ker(\beta^*)\srel{\alpha^*}{\lar}\pi^*\bbL^\vee
\srel{\beta^*}{\lar}
\cal S^*_\ux=\mbox{$\ouset{j=1}{m}{\bigoplus}\eO_{\pi^{-1}(x_j)}(1)$}\to0.
}\label{eq:dual}
\end{m-eqn}
The diagrams below show, respectively, the meaning of equivalent extensions, defining the same class, on the left, and the $\tld G_\ux:= \Aut(\bbL)\times\eO_\ux^\times$-action on $E_\ux$, on the right:
\begin{m-eqn}{
\begin{array}{l|l}
\xymatrix@R=1.5em@C=3.5em{ 
\pi^*\bbL\ar@{^(->}[r]^-{\alpha}\ar@{=}[d] & \cV\ar[d]^-\cong_-{\tld w} \ar@{->>}[r]^-\beta & \cal S_\ux\ar@{=}[d]
\\ 
\pi^*\bbL\ar@{^(->}[r]^-{\alpha'={\tld w}\alpha} & \cV'\ar@{->>}[r]^-{\beta'=\beta{\tld w}^{-1}} & \,\cal S_\ux.}
&
\xymatrix@R=1.5em@C=0.85em{ 
\kern7ex e_\ux: 
&
\pi^*\bbL\ar@{^(->}[rrr]^-{\alpha}\ar[d]_-w 
&&& 
\cV\ar@{=}[d]\ar@{->>}[rrr]^-\beta 
&&& 
\cal S_\ux\ar[d]^-t
\\ 
(w,t)\times e_\ux: 
&
\pi^*\bbL\ar@{^(->}[rrr]^-{\alpha'=\alpha w^{-1}} 
&&& 
\cV\ar@{->>}[rrr]^-{\beta'=t\beta} 
&&& 
\,\cal S_\ux.}
\end{array}
}\label{eq:wse}
\end{m-eqn}
If $\cV'=\cV$ is stable, $\tld w$ is the multiplication by some $c\in\Bbbk^*$. Although $\Bbbk^*_{diag}\subset\tld G_\ux$ acts trivially on extension classes, it acts by rescaling on the vector bundles themselves.

\begin{m-lemma}\label{lm:Ex}
Suppose that $e_\ux\in E_\ux$ is generic that is, it determines a stable bundle $\cV$.

\nit{\rm(i)} The stabilizer of $e_\ux$ is the diagonally embedded $\Bbbk^*_{diag}\subset\tld G_\ux$, so $G_\ux:=\tld G_\ux/\Bbbk^*_{diag}$ acts on $E_\ux$.

\nit{\rm(ii)} Suppose $e_\ux, e'_\ux\in E_\ux$ determine isomorphic $\cV, \cV'$. Then they are in the same $G_\ux$-orbit. 
\end{m-lemma}

\begin{m-proof}
(i) The automorphisms of $\cV$ are multiplications by $c\in\Bbbk^*$. The conclusion follows from the second diagram~\eqref{eq:wse}.

\nit(ii) The isomorphism $\cV\srel{\tld w}{\to}\cV'$ induces $\pi_*\tld w:\bbL=\pi_*\cV\to\pi_*\cV'=\bbL$ that is, $w\in\Aut(\bbL)$. At quotient level, we obtain $t:\cal S_\ux=\cV/\pi^*\bbL\to\cV'/\pi^*\bbL=\cal S_\ux$. 
\end{m-proof}

\subsubsection{The relative case} 
To describe the situation for variable $\ux\in\cA$, we consider the diagram: 
$$
\xymatrix@C=3em@R=1.5em{
\cal X\ar[r]^-{\si'}\ar[d]_-\xi
&
\cal Z\;\ar[d]_-{\zeta}\ar@{^(->}[r]\ar[rd]|(.55){\;{\rm pr}_{\mbb P^1}}
&
{\rm Hilb}^m_{\mbb P^1}\times\mbb P^1\ar[d]
&
Y_\ell\ar[ld]^-{\pi}
\\ 
{(\mbb P^1)}^m\ar[r]^-\si
&
{(\mbb P^1)}^m/\mfrak S_m=\mbb P^m\cong{\rm Hilb}^m_{\mbb P^1}&\mbb P^1&
}
$$
Here $\cal Z$ is the universal family on ${\rm Hilb}^m_{\mbb P^1}$, $\cal X:=\cal Z\times_{\mbb P^m}(\mbb P^1)^m$, and $\mfrak S_m$ are the permutations of $m$ elements. 
In this setting,  $E_\ux$ is the stalk at $\ux$ of the locally free sheaf of rank $2mr$:
$$
E:=
\zeta_*\Bigl({\rm pr}_{\mbb P^1}^*\bigl(\,\cHom(\mbb L^\vee, \pi_*\eO_{\pi}(1))\,\bigr)\Bigr)
=
\zeta_*\Big({\rm pr}_{\mbb P^1}^*\big(\,
\mbb L\otimes(\eO_{\mbb P^1}\oplus\eO_{\mbb P^1}(\ell))
\,\big)\Big).
$$
\begin{m-remark}
For conciseness, we identify $E$ with the linear fibre bundle ${\rm Spec}({\rm Sym}^\bullet_{\eO} E^\vee)$ over $\cA$. 
We are eventually interested in birational properties of $E$, so \textit{$\cA\subset{\rm Hilb}^m_{\mbb P^1}$ is allowed to shrink further}. In the sequel, we denote $E^\cA:=E\rst_\cA$ and $\cU:=\si^{-1}(\cA)\subset{(\mbb P^1)}^m$.
\end{m-remark}

We trivialize $\eO_{\mbb P^1}(-a), \eO_{\mbb P^1}(-1), \eO_{\mbb P^1}(\ell)$ appearing in $\bbL$ and $E$ over $\mbb A^1=\mbb P^1\setminus\{\infty\}$, so their pull-back by ${\rm pr}_{\mbb P^1}$ to $\cal Z_\cA$ is identified with $\eO_{\cal Z_\cA}$:  
\begin{m-eqn}{
\begin{array}{l|l}
\mbb L\rst_{\mbb A^1}\cong
\eO_{\mbb A^1}^{\oplus r-\rho}\oplus
\eO_{\mbb A^1}^{\oplus \rho}=\eO_{\mbb A^1}^{\oplus r},
&\textrm{(View as $r=\!(r-\rho)\!+\!\rho$-column vector.)}
\\[1ex]
E\rst_\cA\cong \zeta_*\bigl(
\eO_{\cal Z_\cA}^{\oplus r}\otimes\eO_{\cal Z_\cA}\oplus
\eO_{\cal Z_\cA}^{\oplus r}\otimes\eO_{\cal Z_\cA}
\bigr) 
&\textrm{(View summands as $r\times m$-matrices.)}
\\[.5ex]
=(\zeta_*\eO_{\cal Z_\cA})^{\oplus r}\oplus (\zeta_*\eO_{\cal Z_\cA})^{\oplus r}
\!=E_\lft\oplus E_\rgt
&E_\lft=\cA\times\bbA^{rm}_\lft, E_\rgt=\cA\times\bbA^{rm}_\rgt
\\[.5ex] 
={(\zeta_*\eO_{\cal Z_\cA}\oplus\zeta_*\eO_{\cal Z_\cA})}^{\oplus r}
={(F\oplus F)}^{\oplus r}.
&F:=\zeta_*\eO_{\cal Z_\cA}=\cA\times\bbA^m.
\end{array}
}\label{eq:Ezo}
\end{m-eqn}
The diagram below describes the situation algebraically: 
$$
\xymatrix@C=1.5em@R=1.5em{
\displaystyle\frac{\Bbbk[x_1,\ldots,x_m][z]}{\langle (z-x_1)\cdot\ldots\cdot(z-x_m)\rangle}
&\displaystyle
\frac{\Bbbk[s_1,\ldots,s_m][z]}{\langle z^m-s_1z^{m-1}+s_2z^{m-2}-\ldots\rangle}
\ar[l]
&
\Bbbk[s_1,\ldots,s_m][z]\ar[l]
\\ 
\Bbbk[x_1,\ldots,x_m]\ar[u]
&\Bbbk[s_1,\ldots,s_m]\ar[u]\ar[l]_-{\text{inclusion}}
&\Bbbk[z],\ar[ul]\ar[u]
}$$
where $s_1=x_1+\ldots+x_m,\ldots,s_m=x_1\cdot\ldots\cdot x_m$ are the symmetric polynomials.  We have: 
\begin{m-eqn}{
\displaystyle 
F=\! 
\frac{\Bbbk[s_1,\ldots,s_m][z]}{\langle z^m-s_1z^{m-1}+s_2z^{m-2}-\ldots\rangle}
\!=
\Bbbk[s_1,\ldots,s_m]\oplus\ldots\oplus
\hat z^{m-1}\!\cdot\Bbbk[s_1,\ldots,s_m]
\cong\eO_{\cA}^{\oplus m},
}\label{eq:zo}
\end{m-eqn}
where $\hat z$ the image of $z$. (The multiplicative structure on $\eO_\cA^m$ is induced by the quotient.) Thus points of $E$ are represented by \textit{pairs} of $r\times m$ block-matrices, with entries in $\Bbbk[s_1,\dots,s_m]$:
\begin{m-eqn}{
\begin{array}{ll}
\scalebox{.9}{$
e=
\left[
\begin{array}{l|l|l|c}
[\text{I}]_{(r-\rho)\times (r-\rho)}
&[\text{III}]_{(r-\rho)\times \rho}
&[\text{V}]_{(r-\rho)\times \rho}&\cdots
\\[1ex] 
 \hline 
\lower1ex\hbox{$[\text{II}]_{\rho\times (r-\rho)}$}
&\lower1ex\hbox{$[\text{IV}]_{\rho\times \rho}$}
&\lower1ex\hbox{$[\text{VI}]_{\rho\times \rho}$}&
\lower1ex\hbox{$\cdots$}  
\end{array}
\right].
$}
&
\begin{array}{l}
\text{The columns of $e$ are:}
\\ 
col_j(e)=\genfrac{[}{]}{0pt}{1}{u_j}{v_j}, 0\les j\les m-1.
\end{array}
\end{array}
}\label{eq:xi}
\end{m-eqn}

\subsection{Groups and slices} 

There are two actions preserving the projection $E\to{\rm Hilb}^m(\mbb P^1)$.

\subsubsection{First symmetry}  $\Aut(\bbL)=\Aut(\bbL^\vee)\subset\End(\bbL)$. It's a linear algebraic group of dimension $r^2$, a representation is obtained by making it act on $\Gamma(\bbL^\vee)\cong\Bbbk^{r-\rho}\oplus\Bbbk^{\rho}\oplus\Bbbk^{\rho}$: 
$$
w=\biggl[
\begin{array}{cc}A&H_0+zH_1\\ 0&B\end{array}
\biggr]
\longmapsto 
\Biggl[
\begin{array}{ccc}A&H_0&H_1\\ 0&B&0\\ 0&0&B
\end{array}
\Biggr],
\quad 
\begin{array}{c}
A\in\Gl(r-\rho), B\in\Gl(\rho),\\ 
H_0, H_1\in\Hom(\Bbbk^{r-\rho},\Bbbk^{\rho})
\end{array}.
$$
Subsequently, we will encounter the following subgroups  of $\Aut(\bbL)$:
$$
U_0^*:=\biggl\{\biggl[
\begin{array}{cc}A&H_0\\0&B\end{array}
\biggr]
\biggr\},
\quad 
U_1^*:=\biggl\{\biggl[
\begin{array}{cc}\bone&zH_1\\0&\bone\end{array}
\biggr]\biggr\}.
$$ 
Let $P\Aut(\bbL)$ stand for $\Aut(\bbL)/\Bbbk^*$, the projective automorphisms. For $e\in E$ as in~\eqref{eq:xi}, $w\times e$ has the same block-form~\eqref{eq:xi}, with the following entries (see~\cite[\S3.3.1]{hlc-sheaves}): 
\begin{m-eqn}{
\begin{array}{rllrl}
w\times{\rm[I]}&=A{\rm[I]}+H_0{\rm[II]}+H_1{\rm[II']},
&&
w\times{\rm[III]}&=A{\rm[III]}+H_0{\rm[IV]}+H_1{\rm[IV']},
\\ 
w\times{\rm[IV]}&=B{\rm[IV]},
&&
w\times{\rm[V]}&=A{\rm[V]}+H_0{\rm[VI]}+H_1{\rm[VI']}.
\end{array}
}\label{eq:wxe}
\end{m-eqn}
The entries ${\rm[II']}, {\rm[IV']}, {\rm[VI']}$ are universal linear combinations of the columns of $w$, see \lcit 
 
\begin{m-lemma}\label{lm:u1-slice}
The $U_1^*$-orbit of the generic $e\in E_\lft$ intersects in a unique point the subspace $\Xi'_U$ defined by the condition $\{{\rm[III]}=0\}$. 
Thus $\Xi'_U\subset E_\lft$ is a slice for the $U_1^*$-action; it's a vector bundle over an open subset $\cA$ of ${\rm Hilb}^m(\mbb P^1)$.
\begin{flushleft}
Recall~\cite[\S3.3.1]{hlc-sheaves}:
\begin{minipage}[t]{0.75\textwidth}
 $\Xi'_A:=\bigl\{{\rm[I]}=c\cdot\bone_{r-\rho}, {\rm[III]}=c\cdot\bone_{\rho}, {\rm[III]=[V]}=0\mid c\in\Bbbk^*\bigr\}\subset E_\lft$\\ 
is a $(\Aut(\bbL),\Bbbk^*_{diag})$-slice, every generic $\Aut(\bbL)$-orbit intersects $\Xi'_A$ along a $\Bbbk^*$-orbit. 
It is a $(U_0^*,\Bbbk^*_{diag})$-slice for the $U_0^*$-action on $\Xi'_U$.
\end{minipage}
\end{flushleft}
\end{m-lemma}

\begin{m-proof}
The explicit form of ${\rm[IV']}$ is (see \lcit): 
$$
{\rm[IV']}
={\rm col}\bigl[v_{r-\rho-1},\dots,v_{r-2}\bigr]
+{\rm col}\bigl[
(-1)^{m-r+\rho-1}s_{m-r+\rho}\cdot v_{m-1},\dots,(-1)^{m-r}s_{m-r+1}\cdot v_{m-1}\bigr].
$$
For generic $e$, it's invertible. (\textit{E.g.} let the first term be the identity and $v_{m-1}=0$.) The matrix $w\in U_1^*$ which cancels the ${\rm[III]}$-component of $e$ has $H_1=-{\rm[III]}\cdot{\rm[IV']}^{-1}$.  
\end{m-proof}


\subsubsection{Second symmetry} 
$\bbT:=(\zeta_*\eO_{\cal Z_\cA})^\times$ is the group scheme of invertible elements in the sheaf of algebras $\zeta_*\eO_{\cal Z_\cA}$. The $\bbT$-action doesnt't preserve the direct summands~\eqref{eq:zo}; \textit{e.g.} $z$ is invertible, it maps the $\hat z^j$- into the $\hat z^{j+1}$-component (except for $j=n-1$, which becomes a linear combination of the previous ones). In matrix terms, $\bbT$ acts on the columns of~\eqref{eq:xi}, but doesn't preserve them individually. Note however, that $\Bbbk[s_1,\ldots,s_m]$ acts componentwise, so $\bbT$ contains a diagonally embedded copy of $\eO_\cA^\times$ (denoted abusively $\Bbbk^*_{diag}$ in the sequel). 

We give now a different description of the $\bbT$-action, over $({\mbb P^1})^m$ rather than ${\rm Hilb}^m_{\mbb P^1}$. Note that $\bbT$ acts diagonally of $F^2=F_\lft\oplus F_\rgt$ --see~\eqref{eq:Ezo}-- and the action on $E=(F^2)^{\oplus r}$ is obtained by repeating it $r$ times. Thus we need to describe the $\si^*\bbT$-action on $\si^*F=\xi_*\eO_{\cal X}$. 

Since  $z^m-s_1z^{m-1}+\dots=(z-x_1)\cdot\dots\cdot(z-x_m)$, we deduce: 
$$
\begin{array}{rl}
\si^*\bbT=(\xi_*\eO_{\cal X})^\times\cong
\big(\eO_{(\mbb P^1)^m}^\times\big)^m,
&
\\[1ex] 
\displaystyle 
\frac{\Bbbk[x_1,\ldots,x_m][z]}{\langle (z-x_1)\cdot\ldots\cdot(z-x_m)\rangle}
&
\cong\Bbbk[x_1,\ldots,x_m]^{\oplus m}\quad\text{is a ring isomorphism,}
\\[1ex]
&
=
{\rm pr}_{\mbb P^1}^*
\big(\eO_{\mbb P^1,x_1}\oplus\ldots\oplus\eO_{\mbb P^1,x_m}\big).
\end{array}
$$
Thus $(t_1,\ldots,t_m)\in(\Bbbk^*)^m$ acts by $t_j$ on the $j$-coordinate of $\si^*F=\mbb A_\cU^m$; it's the action~\eqref{eq:wse}. But now we must also take into account the permutation group $\mfrak S_m$ which interchanges the factors of ${(\mbb P^1)}^m$. For a reason to be clarified (cf. Lemma~\ref{lm:qT-slice}), we consider the affine line $\mbb A^1_\cU\to\cU$ endowed with trivial actions of $\mfrak S_m, \si^*\bbT, \Aut(\bbL)$. 
We add it to $\si^*E_\rgt=(\si^*F_\rgt)^r$, so we obtain the locally trivial linear bundle (vector bundle): 
$$
\tld F:=\si^*E\oplus\mbb A^1_\cU= \si^*E_\lft\oplus \si^*E_\rgt \oplus\mbb A^1_\cU. 
$$
Let $\mu_\Bbbk=\Bbbk^*$ be the multiplicative group, let it act diagonally on $\si^*E_\rgt\oplus\mbb A^1_\cU$. We construct a slice for the $\bbT$-action on $E$ in two stages: 

\smallskip\nit\underbar{Step~1}\quad 
Consider the obvious $\mfrak S_m\times \mu_\Bbbk$-invariant linear subspace in $\tilde F$:
$$
\wtld\Xi'':=\si^*E_\lft\oplus \wtld\Xi''_\rgt,\quad\text{where}\;\; 
\wtld\Xi''_\rgt:=
\scalebox{1.35}{$\Biggl\{\Biggl[$}
\underbrace{
\scalebox{.9}{$
\begin{array}{c|ccc|c|}c&&\ldots&&c\\ *&&\ldots&&*\\ \vdots&&\ldots&&\vdots
\end{array}
$}}_{\Xi''_\rgt\subset\; \si^*E_\rgt}
\scalebox{.9}{$
\begin{array}{c}c\\ 0\\ \vdots
\end{array}
$}
\scalebox{1.35}{$\Biggr]\;\Biggl|\;$}
c\in\Bbbk^*
\scalebox{1.35}{$\Biggr\}$}.
$$
It's a (locally trivial) vector bundle over $\cU$, of relative codimension $m$. 
Note that  $\mu_\Bbbk\!\times\!\si^*\bbT$ intersects $\wtld\Xi''$ along a $\mu_\Bbbk$-orbit. For this reason, we think of $\wtld\Xi''$ as a $\si^*\bbT\!=\!(\mu_\Bbbk\times \si^*\bbT)/{\mu_\Bbbk}$-slice; it's basically the same as setting $c=1$ in $\Xi''_\rgt$ above.

\smallskip\nit\underbar{Step~2}\quad 
Descend to ${\rm Hilb}^m_{\mbb P^1}$ that is, factor out the $\mfrak S_m$-action. This is not automatic, since $\wtld\Xi''$ isn't $\Aut(\bbL)$-invariant and we wish to keep this action. We consider the diagram: 
$$
\xymatrix@R=1.5em@C=3em{
\tld F=\cU\times\bbA^{rm}_\lft\times(\bbA^{rm}_\rgt\times\bbA^1)
\ar@{-->}[r]^-{({\rm id},\Upsilon)}\ar[d]_-{\rm pr}
&
\cal U\times\mbb A^{rm}_\lft\times\times\bbA^{rm+1}_\emptyset
\ar[d]^-{\rm pr}
\\ 
\si^*E_\lft=\cU\times\bbA^{rm}_\lft
\ar@{=}[r]\ar[d]
&
\cU\times\bbA^{rm}_\lft
\ar[d]
\\ 
\cU 
\ar@{=}[r]
&
\cU 
}
$$
The assumptions of the no-name lemma~\cite{bo+ka,dom,dol} are satisfied: $\mfrak S_m\!\times\!\Aut(\mbb L)$ acts on the fibre bundle $\tld F\to \si^*E_\lft$, the generic stabilizer on $\si^*E_\lft$ is trivial. 
Thus there is a $\mfrak S_m\!\times\!\Aut(\mbb L)$-invariant open subset $\tld O\subset\cal U\times\mbb A^{rm}_\lft$, and a birational ${\rm pr}$-fibrewise linear map $({\rm id},\Upsilon)$ which is equivariant for the following actions: 
\begin{itemize}[leftmargin=3ex]
\item $\mfrak S_m$ acts trivially on $\bbA^{rm}_\lft$ (since $\si^*E_\lft$ is pulled-back from the quotient);
\item 
$\mfrak S_m$ acts on $\mbb A^{rm}_\rgt$ by permuting the 
$m$ copies of $\mbb A^{r}$;
\item 
$\Aut(\mbb L)$ acts the same way on $\mbb A^{rm}_\rgt$ and $\mbb A^{rm}_\lft$. (It acts diagonally on $\mbb A^{2rm}$.);
\item 
$\mfrak S_m\times\Aut(\mbb L)$ acts trivially on $\mbb A^{rm+1}_{\emptyset}$.
\end{itemize}
The group $\mu_\Bbbk\times \si^*\bbT$ acts both on the fibre and the base of ${\rm pr}$; the stabilizer of the subvariety $\wtld\Xi''\rst_{\tld O}$ is $\mu_\Bbbk$. A dimensional count shows that its $\mu_\Bbbk\times \si^*\bbT$-orbit is open in $\tld F$; equivalently, the orbit of the generic point in $\tld F$ intersects $\wtld\Xi''\rst_{\tld O}$. 

\begin{m-lemma}\label{lm:qT-slice}
The subvariety 
$\wtld\Xi''_\emptyset\!:=({\rm id},\Upsilon)(\tld\Xi''\rst_{\tld O})\subset
\tld O\!\times\!\mbb A^{rm+1}_\emptyset\!=\mbb A^{rm+1}_{\tld O}$ 
has the properties: 
\begin{enumerate}[leftmargin=5ex]
\item 
It's invariant under the $\mfrak S_m$-action on $\tld O$ and the fibrewise $\mu_\Bbbk$-action over $\tld O$. 
Also, it's a locally trivial linear fibre bundle (vector bundle) over $\tld O$; 
\item 
$\Xi''_\emptyset:=\wtld\Xi''_\emptyset/\mfrak S_m\subset
\big(\tld O/\mfrak S_m\big)\times\mbb A^{rm}_\emptyset\subset E$ 
is a locally trivial vector bundle on $O\!:=\tld O/\mfrak S_m$, and $\Aut(\bbL)$ acts fibrewise on $O$ over $\cA$.
\item $\mbb P(\Xi''_\emptyset)=\Xi''_\emptyset\invq \mu_\Bbbk$ is a locally trivial projective bundle over $O\subset E_\lft$; 
\item[] $\Xi''_T:=({\rm id},\Upsilon)(\Xi''_\rgt\rst_{c=1})/\mfrak S_m\subset E\cap \mbb P(\Xi''_\emptyset)$ is open; it's a slice for the $\bbT$-action on $E$.
\end{enumerate}
\end{m-lemma}

\begin{m-proof}
(i) The invariance follows from fact that $\Xi''_\rgt$ is so. The linearity of $({\rm id,\Upsilon})$ implies that $\wtld\Xi''_\emptyset$ is an linear space bundle over $\si^*E_\lft$. To prove local triviality, take a point $*\in \tld O$ and let $\tld Z_*\subset\mbb A^{rm+1}_{\tld O,*}$ be a complement subspace of $\wtld\Xi''_{\emptyset,*}$; extend it to $\tld O\times\tld Z_*$, trivially. The composed linear map 
$\wtld\Xi''_\emptyset\to\tld F\to\tld F/(\tld O\times\tld Z_*)$
is an isomorphism at $*\in\tld O$, so it's the same in a neighbourhood. But the right-hand side is, obviously, locally trivial. 
For the last claim, we apply the no-name lemma to the group $\mfrak S_m$ acting on $\wtld\Xi''_\emptyset\to\tld O$. The subvariety $\Xi''_T$ is a $\bbT$-slice because $\tld\Xi''\rst_{\tld O}$ is so.
\end{m-proof}


\subsubsection{Summary} 
The group scheme $\tld G:=\Aut(\bbL)\times\bbT$ acts on $E$. The stabilizer of generic extension class is the diagonally embedded $\Bbbk^*_{diag}$. The latter acts by multiplication on the extension themselves. 

The rational map $\,\Psi:E\dashto\mm_{Y_\ell}(r;m)^\svb\,$ is invariant, for the effective action of the group $G:=\tld G/\Bbbk^*_{diag} =(\Aut(\bbL)\times\bbT)/\Bbbk^*_{diag}$. The moduli space of (generic) stable vector bundles is the quotient $E\invq G$ of a suitable open subset. 

\begin{m-lemma}\label{lm:UT-slice}
The restriction $\Xi_{UT}:=\Xi''_T\rst_{O\cap\Xi'_U}\subset E$ is a slice for the $U_1^*\times\bbT$-action.  It is a rational variety, open subset of a locally trivial fibration over $\cA\subset{\rm Hilb}^m_{\mbb P^1}$. 
\begin{flushleft}
Recall~\cite[\S3.3.4]{hlc-sheaves}:  
\begin{minipage}[t]{0.75\textwidth}
$\Xi_G:=\Xi''_T\rst_{O\cap\Xi'_A}$ is a rational slice for the $G$-action.\\ 
It is also a $U_0^*$-slice for the $U_0^*$-action on $\Xi_{UT}$.
\end{minipage}
\end{flushleft}
\end{m-lemma}

\begin{m-proof}
We combine the lemmata~\ref{lm:u1-slice} and~\ref{lm:qT-slice}. The intersection $O\cap\Xi'_U$ is non-empty because $O\subset E_\lft$ is $\Aut(\bbL)$-invariant and the $\Aut(\bbL)$-orbit of $\Xi'_U$ is open in $E_\lft$. 
\begin{flushleft}
\begin{minipage}[t]{0.5\textwidth}
We observe that $\Xi$ fits into the diagram on\\ 
the right, and this proves the second claim.
\end{minipage}
\quad
$\xymatrix@R=1.5em@C=0.25em{
\Xi_{UT}\ar[d]_-{\genfrac{}{}{0pt}{1}{\rm loc.}{\rm triv.}}^-{\genfrac{}{}{0pt}{1}{\rm affine}{\rm bdl.}}&&
\\ 
O\cap\Xi'_U&\srel{open}{\subset}&\Xi'_U\ar[d]_-{\genfrac{}{}{0pt}{1}{\rm loc.}{\rm triv.}}^-{\genfrac{}{}{0pt}{1}{\rm vect.}{\rm bdl.}}
\\ 
&&\cA
}$\vspace{-1.25em}
\end{flushleft}
\end{m-proof}


\subsection{Quotients} 

Let $l_o=\pi^{-1}(0)$ and $\mm_{Y_\ell}(r;m)_{l_o\Lambda}^{\svb}$ be the moduli space of stable vector bundles on $Y_\ell$ which are framed along $l_o\cup\Lambda$; we call the latter $l_o\Lambda$-frames, for short. 

First, note that $l_o\cap\Lambda$ is a point. Thus, if $\cV$ is trivializable along both  $l_o, \Lambda$ (cf. Theorem(ii) above), there is a one-to-one correspondence between $l_o$- and $l_o\Lambda$-frames. Second, it's not clear a priori that $\mm_{Y_\ell}(r;m)_{l_o\Lambda}^{\svb}$ has the structure of a quasi-projective variety. The forthcoming discussion addresses this matter.  

\begin{m-definition}
Let $\cV=\cV_e$ be determined by~\eqref{eq:mihz}, corresponding to the point $e\in E$.  
\begin{enumerate}[leftmargin=5ex]
\item 
We  say that $\cV\in \mm_{Y_\ell}(r;m)_{l_o\Lambda}^{\svb}$ is generic if it's stable, $\cV\rst_{l_o\cup\Lambda}$ is trivializable, and $R^1\pi_*\cV(-\Lambda)$ consists of $m$ distinct points in $\mbb A^1\setminus\{0\}$. Thus $\det(\alpha)\in\Gamma(\eO_{\mbb P^1}(m))$ vanishes at $m$ distinct points, different of $0\in\mbb P^1$. 
Let  $\mm_{l_o\Lambda}^\cA\subset \mm_{Y_\ell}(r;m)_{l_o\Lambda}^{\svb}$ be the corresponding open subset. There is a natural forgetful morphism 
$\mm_{l_o\Lambda}^\cA\to\mm_{Y_\ell}(r;m)^\svb.$ 
\item 
We say that the extension $e\in E$ is generic if it determines a generic vector bundle; let $E^\cA\subset E$ be the open locus determined by generic extensions. 
\item 
We fix an isomorphism 
$\,\Bbbk^r=\Bbbk^{r-\rho}\oplus\Bbbk^\rho\srel{\cong}{\lar}\eO_{0}^{r-\rho}\oplus\eO_0^\rho=\bbL_0\,$ 
which respects the decomposition. Then $\cV$ automatically inherits the framing 
$$
\cV^\vee\to\cV_{l_o}^\vee\srel{\alpha^*(0)}{-\kern-1ex\lar}\pi^*\bbL_0^\vee=\eO_{l_o}^{\oplus r}\;\;
\text{(equivalently, $\alpha(0):\eO_{l_o}^{\oplus r}\to\cV_{l_o}$)} .
$$ 
\end{enumerate}
\end{m-definition}

The automorphisms of $\bbL$ preserve the subspace $\eO_{\mbb P^1}(-a)^{r-\rho}\subset\bbL$, so the frames determined by $\alpha(0)$, as above, don't exhaust all the possible frames $\Bbbk^r\to\cV_{l_o}$.  (The only exception occurs for $\rho=0$ that is, $m=ar$.) Hence, the morphism $E^\cA\to\mm_{l_o\Lambda}^\cA$ is not dominant; also, $U_1^*$ acts trivially on $l_o$-frames. To compensate for this deficiency, we consider the unipotent group 
$$
U_*:=\biggl\{\biggl[
\begin{array}{ll} \bone_{r-\rho}&0\\ \scalebox{1.35}{$*$}_{\rho\times(r-\rho)}&\bone_\rho\end{array}
\biggr]\biggr\},
$$
and define the morphism 
$
\,\Phi:U_*\times_\Bbbk E^\cA\to \mm_{l_o\Lambda}^\cA,\;
(u,e)\mt ([\cV_e], u\cdot\alpha^*(0)).
$

\begin{m-proposition}
\begin{enumerate}[leftmargin=5ex]
\item 
The morphism $\Phi$ is dominant. 
\item 
The group scheme $U_*\times_\Bbbk G$ acts on $U_*\times E$ and $\Phi$ is equivariant for the action.
\item 
The stabilizer in $U_*\times G$ of $(u,e)\in U_*\times_\Bbbk E^\cA$ is trivial.
\item 
If $(u,e), (u',e')$ determine isomorphic framed bundles, $\Phi(u,e)=\Phi(u',e')$, they belong to the same $U_1^*\times_\Bbbk\bbT$-orbit. The generic $U_1^*\times_\Bbbk\bbT$-stabilizer in $U_*\times_\Bbbk E$ is trivial. 
\\ As quasi-projective variety, $\mm_{l_o\Lambda}^\cA$ is the quotient of $U_*\!\times_\Bbbk\!E^\cA$ by the action of $U_1^*\!\times_\Bbbk\!\bbT$.
\end{enumerate}
\end{m-proposition}

\begin{m-proof}
(i),~(ii) The map $\Psi:E^\cA\to\mm_{Y_\ell}(r;m)^\svb$ is dominant, $G$-invariant. The image of $U_*\times U_0^*\to\Gl(r)$ is open, the action preserves $\Psi$, so $\Phi$ exhausts (almost) all the $l_o$-frames. 

\nit(iii) If $(g,\tld w)\in U_*\times\tld G$ stabilizes $(u,e)$, then $\tld w$ stabilizes $e$, so $\tld w=c\in\Bbbk^*$ (cf.~Lemma~\ref{lm:Ex}). Now impose that $g$ stabilizes the framing: $\alpha(0)\cdot g^{-1}=const\cdot\alpha(0)$; it implies $g=1\in U_*$. 

\nit(iv) The vector bundles $\cV, \cV'$ determined by $e, e'$ are isomorphic, so $e'$ is in the $G$-orbit of $e$; let $e'=(w,t)\times e$, with $w\in\Aut(\bbL), t\in\bbT$. Then we have the diagram: 
$$
\xymatrix@R=1.5em@C=5em{ 
\pi^*\bbL\ar@{^(->}[r]^-{\alpha}\ar[d]_-w 
& 
\cV\ar[d]^-\cong_-{\tld w}\ar@{->>}[r]^-\beta 
& 
\cal S_\ux\ar[d]^-t
\\ 
\pi^*\bbL\ar@{^(->}[r]^-{\alpha'=\tld w\alpha w^{-1}} 
& 
\cV'\ar@{->>}[r]^-{\beta'=t\beta{\tld w}^{-1}} 
& 
\,\cal S_\ux.
}
$$
The frames are isomorphic, too, so there is $c\in\Bbbk^*$ such that we have:
$$
\alpha'(0)\cdot u'^{-1}=c\cdot \tld w(0)\cdot\alpha(0)\cdot u^{-1}  
\Rightarrow\; c\cdot w(0)=u'^{-1}u\in U_*\cap\Aut(\bbL)=\{\bone\}.
$$ 
We conclude $u'=u$, $c\cdot w\in U_1^*$ and $t\in\bbT$ is arbitrary. By replacing $\tld w\mt\tld w_{new}=c\cdot\tld w$, similarly for $\alpha,\beta$, we preserve the extension class $e$  and $w_{new}\in U_1^*$. Finally, the generic stabilizer in $\tld G$ is $\Bbbk^*_{diag}$ which intersects $U_1^*\times\bbT$ trivially.
\end{m-proof}

The following commutative diagram of rational morphisms summarizes the situation:
$$
\xymatrix@R=3em@C=5em{
U_*\times E^\cA 
\ar[rr]^-{\Phi}_-{U_1^*\times\bbT}|-{\;quotient\;}
\ar@{=}[d]
&&
\mm_{l_o\Lambda}^\cA
\ar[d]^-{
\scalebox{.8}{$\begin{array}{l}
{\rm PGl}(r)\,\srel{open}{\supset}\,U_*\times PU_0^*\\ \text{acts on frames} 
\end{array}$}}
\\ 
U_*\times E^\cA 
\ar@/_4ex/[rr]^-{\Psi}_-{U_*\times G}|-{\;quotient\;}
\ar[r]^-{{\rm pr}_{E_\cA}}_-{U_*}
&E^\cA\ar[r]^-{quotient}_-{G}&
\mm_{Y_\ell}(r;m)^\svb
}
$$
The decorations indicate the general fibres. The slices $\Xi_{UT}, \Xi_G$ fit into the same diagram (see the `Recall'-comment in Lemma~\ref{lm:UT-slice}).

\begin{m-theorem}\label{thm:rtl-frame}
\begin{enumerate}[leftmargin=5ex]
\item 
The moduli space $\mm_{Y_\ell}(r;m)_{l_o\Lambda}^{\svb}$ is an irreducible, rational variety. Also, it admits a universal Poincar\'e bundle over an open subset. 
\item 
The moduli space $\mm_{Y_\ell}(r;m)$ is an irreducible, rational variety 
{\rm(}cf.~\cite[Theorem 3.8]{hlc-sheaves}{\rm)}.
\end{enumerate}
\end{m-theorem}
Note that the first part of the result is in agreement with Bartocci~\textit{et al.}~\cite{bbr}. 
 
\begin{m-proof}
The first claim follows from the previous Proposition and Lemma~\ref{lm:UT-slice}: a slice for the $U_1^*\times\bbT$-action is $U_*\times\Xi_{UT}$, which is a rational variety. Now we construct the universal framed bundle. Note that $\tld{\cal Z}:=E\times_{\cal H}(\cal Z\times_{\mbb P^1}Y_\ell)$ is a subvariety of $E\times Y_\ell$, since $\cal Z\subset\cal H\times\mbb P^1$. Let $s$ be the tautological section of $\zeta^* E$ over $E$. The diagram below globalizes~\eqref{eq:dual}:
$$
\scalebox{.9}{$
\xymatrix@C=3em@R=3em{
{(\pi\circ{\rm pr}_{Y_\ell}^{E\times Y_\ell})}^*\, \bbL^\vee\ar[r]\ar@/_3ex/[drr]|-{\;\beta^*\;}
&
{(\pi\circ{\rm pr}_{Y_\ell}^{E\times Y_\ell})}^*\, 
\bbL^\vee\otimes\eO_{\tld{\cal Z}}
\ar[r]^-{\langle\cdot,s\rangle}_-{\text{pairing}}
&
{(\pi\circ{\rm pr}_{Y_\ell}^{E\times Y_\ell})}^*\pi_*\eO_\pi(1)\otimes\eO_{\tld{\cal Z}}
\ar[d]^-{\text{evaluation}}
\\&&
{({\rm pr}_{Y_\ell}^{E\times Y_\ell})}^*\eO_\pi(1)\otimes\eO_{\tld{\cal Z}}.
}
$}
$$
The pairing with $s$ is generically surjective, $(U_1^*\times\bbT)\cdot\Xi_{UT}\subset E$ is open, so the vector bundle  $\cal W:=\Ker(\beta^*)|_{\Xi_{UT}\times Y_\ell}$ fits into the universal exact sequence over $\Xi_{UT}\times Y_\ell$: 
$$
0\to\cal W\srel{\alpha^*}{\lar}
{(\pi\circ{\rm pr}_{Y_\ell}^{E\times Y_\ell})}^*\,\bbL^\vee 
\srel{\beta^*}{\lar}
{({\rm pr}_{Y_\ell}^{E\times Y_\ell})}^*\eO_\pi(1)\otimes\eO_{\tld{\cal Z}}
\to0.
$$
To include frames into the picture, we consider the pull back to $U_*\times\Xi_{UT}\times Y_\ell$. Then $\cV^\vee:={({\rm pr}^{U_*\times\Xi_{UT}\times Y_\ell}_{\Xi_{UT}\times Y_\ell})}^*\cal W$ possesses the following universal framing over 
$U_*\times\Xi_{UT}\times l_o$:\\[1ex] 
$f:\cV^\vee\to\cV^\vee\rst_{U_*\times\Xi_{UT}\times l_o}
\srel{\alpha^*(0)}{-\kern-1ex\lar} 
{\bigl({\rm pr}^{U_*\times\Xi_{UT}\times Y_\ell}_{0}\bigr)}^*\bbL_0^\vee 
=\eO_{U_*\times\Xi_{UT}\times l_o}^{\oplus r},
\;\; f_{(u,e)}:=u\cdot \alpha^*_e(0).$
\end{m-proof}

\begin{m-remark}
The first named author claimed~\cite[Corollary 3.11]{hlc-sheaves} the rationality of the framed moduli space. The argument involves a generic Poincar\'e bundle, induced from the $G=\tld G/\Bbbk^*_{diag}$-slice $\Xi_G$ in $E$. Unfortunately, this is incorrect: $\Bbbk^*_{diag}\subset\tld G$, although acts trivially on extension \textit{classes}, it acts by multiplication on the extensions. Briefly, the generic stabilizer acts non-trivially, descent doesn't apply. In contrast, in our framed situation, the $U_1^*\times_\Bbbk\bbT$-stabilizer in $U_*\times E$ is trivial, so descent is applicable in Theorem~\ref{thm:rtl-frame}. 
\end{m-remark}


\subsection{Application to the plane and quadric}\label{ssct:applic}

We apply the results in the following cases:  
\begin{enumerate}[leftmargin=5ex]
\item the $0^{\rm th}$ Hirzebruch surface $Y_0=\cQ$, which is $\mbb P^1\times\mbb P^1$. 
\item[] 
Here the charge $m=2n$ (the quadric has degree $2$ in $\ppp$), so $\mm_\cQ(r;2n)$ is irreducible, rational, of dimension $4rn-r^2+1$. We immediately deduce the form~\eqref{eq:miq} of the general vector bundle in $\mm_\cQ(r;2n)$.\smallskip
\item 
the $1^{\rm st}$ Hirzebruch surface $Y_1$, which is the blow-up of a plane at a point. 

Here the charge is $m=n$. Any semi-stable vector bundle on $\pp$ admits a (semi-stable) deformation whose restriction to the general line is trivializable (Hirschowitz-lemma, see~\cite{hulek}). Thus there is a line $\lda\subset\pp$ such that the restriction to $\lda$ of the general vector bundle in each (possible) irreducible component of $\mm_\pp(r;n)$ is trivializable. 

Let $Y_1\srel{\Pi}{\to}\pp$ be the blow-up of a point $p\in\lda$. Then $\Pi^*$ and $\Pi_*$ determine birational maps $\mm_{Y_\ell}(r;n)\dashto\mm_\pp(r;n)$. The theorem above yields the irreducibility and rationality of $\mm_\pp(r;n)$; it's dimension is $2nr-r^2+1$. (See also~\cite[Corollary 3.9]{hlc-sheaves}. Hulek~\cite{hulek} obtained the irreducibility in a different way.) 
The general stable vector bundle on $Y_1$ is given by~\eqref{eq:mihz}; since 
$\;\pi^*\eO_{\mbb P^1}(1)=\Pi^*\eO_\pp(1)\otimes\eO_{Y_1}(-\Lambda)$ and $\; \Pi_*\eO_{Y_1}(-\Lambda)=\eI_p\subset\eO_\pp,$ 
we obtain the exact sequence~\eqref{eq:mip}.
 \end{enumerate}
 

\section{Returning to the roots}\label{sct:roots}

In this final section, we investigate the implications of our results to the Yang-Mills instantons, which have been constituting the motivation for investigating their mathematical generalization. 

\subsection{Consequences for physical instantons}\label{ssct:ath}

Our work is related to an issue raised by Atiyah, solved by Donaldson; it is listed in Hartshorne~\cite[Problem~22]{hart-probl}. 
We briefly recall the setup (see~\cite{atiy-2+4,atiyah-ward,dons,hart-inst} for details). 

Consider the field of quaternions acting on itself by left-multiplication: 
$$
\begin{array}{rl}
\mbb H&=\{r_1+r_2\bi+r_3\bj+r_4\bk\mid r_1,\dots, r_4\in\mbb R\}=\mbb R^4\\ 
&=\{ z_1-z_2\bj\mid z_1, z_2\in\mbb C \}=\mbb C^2.
\end{array}
$$
The twistor fibre bundle map $t\!w$ is defined as follows:
$$
\xymatrix@R=1.5em@C=1.5em{
\mbb C^4\setminus\{0\}\ar[r]\ar[d]
&\mbb H^2\setminus\{0\}\ar[d]&&\\ 
\mbb{CP}^3\ar[r]^-{t\!w}&\mbb{HP}^1=\mbb S^4
&[z_1:z_2:z_3:z_4]\ar@{|->}[r]&(z_1-z_2\bj)\cdot(z_3-z_4\bj)^{-1}\in\mbb H\cup\{ \infty \}=\mbb S^4.
} 
$$
Left-multiplication by $\bj$ on $\mbb C^4$ determines the (real) automorphism 
$$
(z_1,z_2,z_3,z_4)\srel{\rho}{\lmt}(\bar z_2, -\bar z_1, \bar z_4, -\bar z_3), 
$$
which descends to $\mbb{CP}^3$. It defines a real structure that is, $\rho$ is conjugate-linear on $\eO_{\mbb{CP}^3}$ and $\rho^2=\bone_{\mbb{CP}^3}$. Note that $\rho$ has neither fixed points nor fixed $2$-planes. Nevertheless, the fibres of $t\!w$ are $\rho$-invariant, they are complex lines ($2$-dimensional spheres) $\mbb{CP}^1\subset\mbb{CP}^3$, with induced real structure; for this reason, one calls them `real lines'. We fix such a real line $\lda$. 

The real-automorphism $\rho$ determines the `dual-conjugate pull-back' map:
$$
\wtld{\rho}^*:\mi_{\mbb{CP}^3}(r;n)_\lda\to\mi_{\mbb{CP}^3}(r;n)_\lda,
\quad\beF\mt{\overline{\rho^*\beF}}^\vee.
$$ 
(The $\lda$-frame is invariant precisely when it's real.) The Penrose transform --cf. Atiyah-Hitchin-Singer \cite[Theorem 5.2]{ahs}-- identifies irreducible, self-dual $SU(r)$-connections on the real $4$-dimensional sphere $\mbb S^4$ with holomorphic vector bundles on $\mbb{CP}^3$, possessing a real structure, which are trivializable on the fibres of the projection $\mbb{CP}^3\to\mbb S^4$. 

More precisely, the moduli space $\iph_{\mbb{CP}^3}(r;n)_\lda$ of $\lda$-framed physical instantons is a $4rn$ real-dimensional sub-manifold of $\mi_{\mbb{CP}^3}(r;n)_\lda$,  it's the $\wtld{\rho}^*$-fixed locus in $\mi_{\mbb{CP}^3}(r;n)_\lda$, and it consists of those $\beF$ which satisfy the properties: 
\begin{itemize}[leftmargin=3ex]
\item 
It is trivializable on the fibres of $t\!w$, which are real lines in $\mbb{CP}^3$;
\item 
There is an isomorphism $\tld j:\beF\to\wtld{\rho}^*\beF$, trivial on real lines. (It's determined up to multiplication by some $c\in\mbb C, |c|=1$.)
\end{itemize}
Let $\cH$ be a $2$-plane containing $\lda$, and $\cD:=\rho(\cH)$. (\textit{E.g.} take $\lda=\{ [z_1:z_2] \}=t\!w^{-1}(\infty)$, $\cH=\{z_4=0\}, \cD=\{z_3=0\}$.)
We consider the restriction map   
$$
\res_\cH :\iph_{\mbb{CP}^3}(r;n)_\lda\to\mm_\cH(r;n)_\lda.
$$
One readily deduces that the map is injective (see Remark~\ref{rk:tpym}), and Atiyah asked whether $\res_\cH$ is a \textit{real} diffeomorphism. Donaldson~\cite{dons} answered this in affirmative, and his proof passes through the Kempf-Ness theory. Furthermore,  Atiyah~\cite[Theorem 1]{atiy-2+4} proved that the real manifold $\iph_{\mbb{CP}^3}(r;n)_\lda$ possesses a natural complex analytic structure, so that $\res_\cH$ is actually bi-holomorphic. 

In the sequel, we show first that our result yields a criterion for recognizing physical (YM) instantons among mathematical ones. Second, we deduce an unexpected relationship between mathematical- and YM-instantons. 

We consider the map, analogous to $\wtld{\rho}^*$ above, and denoted the same:  
$$
\wtld{\rho}^*:\mm_\cH(r;n)_\lda\to\mm_\cD(r;n)_\lda,\quad \cV_\cH\mt\cV_\cD:={\overline{\rho^*\cV_\cH}}^\vee.
$$ 
Let $\Delta_{\wtld{\rho}^*}\subset\mm_\cD(r;n)_\lda\times\mm_\cH(r;n)_\lda$ be the graph; it is the  $\wtld{\rho}^*$-fixed locus in the product. 

\begin{m-lemma}\label{lm:real-inst}
A mathematical instanton $\beF\in\mi_{\mbb{CP}^3}(r;n)_\lda$ is a physical instanton if and only if 
the restriction of $\beF$ to $\cD\cup\cH$ belongs to $\Delta_{\wtld{\rho}^*}$. 
\end{m-lemma}

\begin{m-proof}
The condition is necessary: YM-instantons satisfy $\wtld{\rho}^*\beF\cong\beF$. It's sufficient: suppose $\beF\rst_{\cD\cup\cH}\cong\wtld{\rho}^*\beF\rst_{\cD\cup\cH}$. Since $H^1(\Hom(\beF,\wtld{\rho}^*\beF)(-2))=0$, the isomorphism extends to $\mbb{CP}^3$. 
\end{m-proof}

This observation allows recovering Donaldson's result from a different viewpoint, by embedding it into a `wider context'. 

\begin{m-theorem}\label{thm:ym}
\begin{enumerate}[leftmargin=5ex]
\item 
The map $\res_\cH$ is a diffeomorphism. 
\item 
Let  $\iph_{\mbb{CP}^3}(r;n)_\lda$ be endowed with the complex structure determined by $\res_\cH$, equivalently, with that provided by Atiyah~\cite{atiy-2+4}. Then $\iph_{\mbb{CP}^3}(r;n)_\lda$ is a complex, irreducible, rational quasi-projective variety of complex dimension $2rn$. 
\end{enumerate}
\end{m-theorem}

\begin{m-proof}
(i) Note that $\res_\cH$ has Zariski-open image in $\mm_\pp(r;n)_\lda$, dense in the analytic topology. 
The morphism $\Theta_{\cD\cH,\lda}$ is open --it's differential is isomorphism-- so its image is open in $\mm_\cD(r;n)_\lda\times\mm_\cH(r;n)_\lda$. 
Thus $V:={\rm pr} (\Delta_{\wtld{\rho}^*}\cap{\rm Image}(\Theta_{\cD\cH,\lda}))\subset\mm_\cH(r;n)_\lda$ is Zariski-open.

For any $\cV_\cH\in V$, the pair $(\cV_\cD:=\wtld{\rho}^*\cV_\cH,\cV_\cH)$ belongs to ${\rm Image}(\Theta_{\cD\cH,\lda})$, so there is a mathematical instanton  $\beG$ whose restrictions to $\cH$ and $\cD$ are $\cV_\cH$ and $\cV_\cD$, respectively. Thus $\beG\rst_{\cD\cup\cH}=\wtld{\rho}^*\beG\rst_{\cD\cup\cH}$, and $\beG$ is a physical instanton. 

For the surjectivity of $\res_\cH$, it suffices to prove the following:

\nit\underline{Claim}\quad For any $\cV\in\mm_\cH(r;n)$ trivializable along $\lda$, there is $\beF\in\mi_\ppp(r;n)$ which is $\wtld{\rho}^*$-invariant and $\beF\rst_\cH=\cV$. 

This follows from a general GIT-argument applied to the quiver below, where the vertical dots indicate $a=h^0(\eO_\ppp(1))=4$ arrows. We indicate the main points and skip details.
$$
Q_a: 
\xymatrix@R=1.5em@C=7em{
\bullet
\ar@/^2.5ex/[r]\ar@/_2.5ex/[r]
\ar@{}[r]|-{\raise1.5ex\hbox{$a$}\;\;\raise1ex\hbox{$\vdots$}\;\;\raise1.5ex\hbox{arrows}}
&\bullet 
\ar@/^2.5ex/[r]\ar@/_2.5ex/[r]
\ar@{}[r]|-{\raise1.5ex\hbox{$a$}\;\;\raise1ex\hbox{$\vdots$}\;\;\raise1.5ex\hbox{arrows}}
&\bullet 
}.
$$ 
We consider the dimension vector $(n,r+2n,n)$ and the corresponding representation space: 
$$
\cal R_a= 
\Hom(\mbb C^{n}_\lft,\mbb C^{r+2n})^{\oplus\, a}
\oplus \Hom(\mbb C^{r+2n},\mbb C^{n}_\rgt)^{\oplus\, a}.
$$
The group $\tld A=\Gl(n)\times\Gl(r+2n)\times\Gl(n)$ acts on it. 
The elements of $\cal R_a$ are pairs $\LL=(L_\lft, L_\rgt)$; for $\uz\in\mbb C^4$, let $L_\rgt(\uz):=L_\rgt^{(1)}z_1+\dots+L_\rgt^{(4)}z_4$ and similarly $L_\lft(\uz)$.  

The ADHM and Barth-Hulek construction~\cite{adhm,brth+hulk} imply the following facts: 
\begin{itemize}[leftmargin=3ex]
\item 
$\mi_{\mbb{CP}^3}(r;n)$ is the quotient of $Mond_{\mbb{CP}^3}(r;n)$ by the $\tld A$-action. (See notation in Section~\ref{sct:monad}.) 
The monad which determines an instanton $\beF$ has the property: 
\begin{m-eqn}{
\mbb C^n_\rgt\cong H^1(\beF(-1)), 
\mbb C^n_\lft\cong H^1(\beF^\vee(-1))^\vee
\Rightarrow\;(\mbb C^n_\lft)^*:=\overline{(\mbb C^n_\lft)^\vee}
\cong H^1(\overline{\beF^\vee}(-1)).
}\label{eq:h1}
\end{m-eqn}
\item 
The moduli space $\iph_{\mbb{CP}^3}(r;n)$ is obtained as follows: 
\begin{itemize}[leftmargin=3ex]
\item One identifies 
$\mbb C^{r+2n}\cong(\mbb C^{r+2n})^*$ using a Hermitian structure. 
\item  
The pair $(L_\lft, L_\rgt)$ defines a Yang-Mills instanton precisely when the $\mbb C$-linear map below is surjective and it commutes with the action of the quaternions, for all $\uz\in\mbb C^4$:
$$
\LL(\uz):={L_\lft(\uz)}^*\oplus L_\rgt(\uz):\mbb C^{r+2n}\to
(\mbb C^n_\lft)^*\oplus\mbb C^n_\rgt. 
$$ 
Since $\mbb H$ is generated over $\mbb R$ by $\bi, \bj, \bk=\bi\bj$ and due to the isomorphisms~\eqref{eq:h1} above, this property is equivalent to the following (where $(\;)^*$ stands for the adjoint): 
\begin{m-eqn}{
\begin{array}{ll}
&\LL(\bi\uz)=\bi\LL(\uz),\;\;\LL(\bj\uz)=\bj\LL(\uz),\;\;\forall\,\uz\in\mbb C^4, 
\\ 
\Leftrightarrow&
(L_\lft^{(1)})^*=-L_\rgt^{(2)},
\;\;(L_\lft^{(2)})^*=L_\rgt^{(1)},
\;\;(L_\lft^{(3)})^*=-L_\rgt^{(4)},
\;\;(L_\lft^{(4)})^*=L_\rgt^{(3)}. 
\end{array}
}\label{eq:j-invar}
\end{m-eqn}
\item 
Note that $Mond_{\mbb{CP}^3}(r;n)\subset Cplx_{\mbb{CP}^3}(r;n)\subset \cal R_a$.
The identities above make sense for complexes and for elements of $\cal R_a$; let $Cplx_{\mbb{CP}^3}(r;n)^{\mbb H}$ be the corresponding locus. 

The quiver $Q_3$ has no cycles, so the invariant quotient $\cal R_a\invq\tld A$ is projective (for any linearization). Therefore  $Cplx_{\mbb{CP}^3}(r;n)\invq\tld A$ is projective too, and contains the closed,  thus compact, $\bj$-invariant locus defined by $Cplx_{\mbb{CP}^3}(r;n)^{\mbb H}$. 
\end{itemize}
\end{itemize}
Back to the claim above: by compactness, there is $\LL=(L_\lft,L_\rgt)\in Cplx_{\mbb{CP}^3}(r;n)^{\mbb H}$ whose restriction to $\cH$ --this amounts to forgetting the $z_4$-components-- is a complex whose cohomology is $\cV$. (We avoid the Kempf-Ness theory.) Since $\cV$ is a (locally free) vector bundle, the restricted complex is actually a monad, $\LL(\uz)\rst_{z_4=0}$ is surjective. The $\mbb{H}$-invariance property~\eqref{eq:j-invar} implies that $\LL(\uz)$ is surjective for all $[\uz]\in\mbb{CP}^3$. (This argument is taken from Donaldson~\cite[pp. 457]{dons}.)

\nit(ii) 
The map $\res_\cH$ is birational and $\mm_{\mbb{CP}^2}(r;n)_\lda$ is irreducible, rational (cf. Theorem~\ref{thm:rtl-frame}).
\end{m-proof}

We continue with a further --novel, to our knowledge-- consequence of Theorem~\ref{thm:birtl}. 

\begin{m-theorem} \label{thm:n2n}
Let $\lda, \lda', \cD, \cH, \cQ$ as in Section~\ref{ssct:res}. There are (algebraic) open immersions: 
$$
\begin{array}{l}
\mi^n_{n,n}:\mi_{\mbb{CP}^3}(r;n)_{\lda}
\to
\iph_{\mbb{CP}^3}(r;n)_\lda\times\iph_\ppp(r;n)_\lda,
\\[1ex] 
\mi^n_{2n}:\mi_{\mbb{CP}^3}(r;n)_{\lda\cup\lda'}
\to
\iph_{\mbb{CP}^3}(r;2n)_{\rm line}.
\end{array}
$$ 
\end{m-theorem}
The second morphism seems especially interesting: leaving frames aside, it says that mathematical instantons of change $n$ are the same as Yang-Mills instantons of charge $2n$.  This matter is not at all obvious in linear algebraic terms, at the monad/ADHM-construction level. The proof shows also that $\mi^n_{2n}$ commutes with the `functors' $\oplus, \otimes$ in Theorem~\ref{thm:categ}, so we have the commutative diagram: 
\begin{m-eqn}{
\xymatrix@R=3em@C=5em{
\mi_{\mbb{CP}^3}(r';n')_{\lda\cup\lda'}\times \mi_{\mbb{CP}^3}(r'';n'')_{\lda\cup\lda'}
\ar[d]_-\otimes 
\ar[r]^-{(\mi^{n'}_{2n'},\mi^{n''}_{2n''})}
&
\iph_{\mbb{CP}^3}(r',2n')_{\rm line}\times \iph_{\mbb{CP}^3}(r'',2n'')_{\rm line}
\ar[d]^-\otimes 
\\ 
\mi_{\mbb{CP}^3}(r'r'';r'n''+r''n')_{\lda\cup\lda'}
\ar[r]^-{\mi^{r'n''+r''n'}_{2(r'n''+r''n')}} 
&
\iph_{\mbb{CP}^3}(r'r'',2(r'n''+r''n'))_{\rm line}.
}}\label{eq:categYM}
\end{m-eqn}
(The Penrose transform implies that the tensor product preserves YM-instantons, so the rightmost vertical arrow is well-defined.)

\begin{m-proof}
We define $\mi^n_{n,n}$ as the following composition:
$$
\mi_{\mbb{CP}^3}(r;n)_\lda\srel{\Theta_{\cD\cH,\lda}}{\arr{7ex}}
\mm_\cD(r;n)_\lda\times\mm_\cH(r;n)_\lda
\srel{(\res_\cD^{-1},\res_\cH^{-1})}{\arr{10ex}}
\iph_{\mbb{CP}^3}(r;n)_\lda\times\iph_{\mbb{CP}^3}(r;n)_\lda.
$$
The first arrow is an open immersion and the second is bi-regular.

The morphism $\mi^n_{2n}$ is defined as follows: 
$$
\mi_{\mbb{CP}^3}(r;n)_{\lda\cup\lda'}
\srel{\Theta_{\cQ,\lda\cup\lda'}}{\arr{7ex}}
\mm_\cQ(r;2n)_{\lda\cup\lda'}
\ouset{\cong}{\alpha}{\lar}
\mm_{\mbb{CP}^2}(r;2n)_{\rm line}
\ouset{\cong}{\res^{-1}}{\lar}
\iph_{\mbb{CP}^3}(r;2n)_{\rm line}.
$$
The isomorphism $\alpha$ is described by Atiyah~\cite[eq. (3.6)]{atiy-2+4},  it's determined by the diagram: 
$$
\xymatrix@R=1.5em@C=3em{
&\tld\cQ
\ar[ld]_(.7){\scalebox{.75}{$\begin{array}{r} 
\text{blow-up}\\ \lda\cap\lda'\end{array}$}\kern-2ex}
\ar[rd]^(.7){\kern-1ex\scalebox{.75}{$\begin{array}{l} 
\text{blow}\\ \text{down}\;\tld\lda, \tld\lda'\end{array}$}}&
\\ 
Q&&{\mbb{CP}^2}.
}
$$
The line denoted `line' above is the image in ${\mbb{CP}^2}$ of the exceptional divisor of $\tld\cQ$, and $\tld\lda, \tld\lda'\subset\tld\cQ$ are the proper transforms of $\lda, \lda'$, respectively. 
\end{m-proof}


\subsection{Final thoughts}

We conclude with a few speculative remarks. Instantons were introduced by physicists ('t~\!Hooft, Polyakov, etc.) motivated by physical phenomena. Thus one naturally wonders whether the brief statement 
\\[1ex]\centerline{`\textit{The tensor product of two mathematical instantons is still a mathematical instanton.}'}\\[1ex] 
has any relevance in physics. 

The authors were very pleased to find that tensor product representations of products of compact groups --in our situation for $SU(r)\times SU(r')$-- has been indeed investigated in the particle physics literature~\cite{cgt,dri-man,jack,be-st,tem}, where the resulting objects are apparently known as `multi-instantons'. The consideration seems to be restricted only to those instantons which originate from the $4$-dimensional sphere through the ADHM-construction. On the physics side, the practical reason for investigating tensor products relies in the computation of their Green functions, which are used for estimating instanton effects in quantum chromodynamics (cf.~\cite[pp. 94]{cgt}). This circle of ideas is beyond the authors' expertise but, keeping in mind the categorical behaviour~\eqref{eq:categYM}, we believe that it's worth mentioning the matter. 

Often, one is concerned with tensor products of vector bundles possessing additional structure (mostly orthogonal or symplectic); in other words, besides $SU(r)$, one is interested in orthogonal or symplectic vector bundles. (For the passage from complex groups to their compact forms, one applies the Kobayashi-Hitchin correspondence.) The statement above remains valid, because the tensor product breaks into irreducible components and the instanton condition --that is, $H^1(\beF\otimes\beG(-2))=0$-- holds for all the direct summands.


\end{document}